\documentclass[11pt]{amsart}

\usepackage{epigamath}

\usepackage[english]{babel}

\usepackage{amscd, amssymb, amsmath, mathrsfs, amsthm}
\usepackage[utf8]{inputenc}
\usepackage{tikz-cd}
\usetikzlibrary{decorations.pathreplacing}
\usepackage{booktabs}
\usepackage{color}

\usepackage{multicol}

\usepackage{microtype}

\theoremstyle{plain}
\newtheorem{theorem}{Theorem}[section]
\newtheorem*{maintheorem}{Theorem}
\newtheorem{lemma}[theorem]{Lemma}
\newtheorem{proposition}[theorem]{Proposition}
\newtheorem{corollary}[theorem]{Corollary}

\theoremstyle{definition}

\theoremstyle{remark}

\DeclareFontFamily{U}{wncy}{}
\DeclareFontShape{U}{wncy}{m}{n}{<->wncyr10}{}
\DeclareSymbolFont{mcy}{U}{wncy}{m}{n}
\DeclareMathSymbol{\Sh}{\mathord}{mcy}{"58}

\newcommand\mylabel[1]{\label{#1}\marginpar{\vspace{-1ex}\medskip\medskip\footnotesize \tt #1}}
\renewcommand\mylabel[1]{\label{#1}}

\newcommand{\bC}{\bar{C}}

\newcommand{\ZZ}	{\mathbb{Z}}
\newcommand{\QQ}	{\mathbb{Q}}

\newcommand{\PP}	{\mathbb{P}}
\renewcommand{\AA}	{\mathbb{A}}
\newcommand{\GG}	{\mathbb{G}}

\newcommand  {\shA}     {\mathscr{A}}

\newcommand  {\shE}     {\mathscr{E}}
\newcommand  {\shF}     {\mathscr{F}}
\newcommand  {\shG}     {\mathscr{G}}

\newcommand  {\shM}     {\mathscr{M}}

\newcommand  {\shN}     {\mathscr{N}}
\newcommand  {\shL}     {\mathscr{L}}

\newcommand{\lieg}{\mathfrak{g}}
\newcommand{\lieh}{\mathfrak{h}}

\newcommand  {\aff}     {{\text{\rm aff}}}

\newcommand  {\alg}     {{\operatorname{alg}}}

\newcommand  {\APic}    {\operatorname{APic}}

\newcommand  {\can}     {{\rm \text{can}}}

\newcommand  {\Cl}      {\operatorname{Cl}}

\newcommand  {\Cokernel}{\operatorname{Coker}}

\newcommand  {\Div}     {\operatorname{Div}}

\newcommand	{\Exc}	{\operatorname{Exc}}
\newcommand  {\End}     {\operatorname{End}}
\newcommand  {\Ext}     {\operatorname{Ext}}

\newcommand  {\Gal}     {\operatorname{Gal}}

\newcommand  {\Hom}     {\operatorname{Hom}}

\newcommand  {\id}      {{\operatorname{id}}}

\newcommand	{\ind}	{\operatorname{ind}}

\newcommand  {\II}  	{\text{\rm II}}

\newcommand  {\Kernel } {\operatorname{Ker}}

\newcommand  {\Lie}     {\operatorname{Lie}}

\newcommand  {\lra}     {\longrightarrow}

\newcommand  {\Mat}     {\operatorname{Mat}}

\newcommand  {\maxid}   {\mathfrak{m}}

\newcommand  {\mult}    {{\text{\rm mult}}}

\newcommand  {\Num}     {\operatorname{Num}}

\newcommand  {\NS}      {\operatorname{NS}}

\renewcommand{\O}       {\mathscr{O}}

\newcommand  {\ord}     {\operatorname{ord}}

\newcommand  {\pdeg}    {\operatorname{pdeg}}

\newcommand  {\Pic}     {\operatorname{Pic}}

\newcommand  {\pr}      {\operatorname{pr}}
\newcommand  {\Proj}    {\operatorname{Proj}}

\newcommand  {\quadand} {\quad\text{and}\quad}

\newcommand  {\ra}      {\rightarrow}

\newcommand  {\rank}    {\operatorname{rank}}
\newcommand  {\red}     {{\operatorname{red}}}
\newcommand  {\Reg}     {\operatorname{Reg}}

\newcommand  {\sep}     {{\operatorname{sep}}}

\newcommand  {\Sing}    {\operatorname{Sing}}

\newcommand  {\Spec}    {\operatorname{Spec}}

\newcommand  {\Supp}    {\operatorname{Supp}}
\newcommand  {\Sym}     {\operatorname{Sym}}

\EpigaVolumeYear{4}{2020}
\EpigaArticleNr{11}
\ReceivedOn{February 24, 2020}
\InFinalFormOn{May 25, 2020}
\AcceptedOn{July 3, 2020}

\title{The maximal unipotent finite quotient, unusual torsion in Fano threefolds, and exceptional Enriques surfaces}
\titlemark{Unusual Torsion in Fano varieties}

\author{Andrea Fanelli}
\address{Institut de mathématiques de Bordeaux (IMB), CNRS, 
Université de Bordeaux,
33405 Talence cedex, France}
\email{andrea.fanelli.1@u-bordeaux.fr}

\author{Stefan Schr\"oer}
\address{Mathematisches Institut, Heinrich-Heine-Universit\"at, 40204 D\"usseldorf, Germany}
\email{schroeer@math.uni-duesseldorf.de}

\authormark{A. Fanelli and S. Schr\"oer}

\AbstractInEnglish{We introduce and study the  maximal unipotent finite quotient for  algebraic   group schemes in  positive characteristics. Applied to  Picard schemes, this quotient encodes  unusual torsion.  We   construct   integral  Fano threefolds where such unusual torsion  actually appears.  The existence of such threefolds is   surprising, because the unusual torsion vanishes for del Pezzo surfaces. Our construction relies on the theory of exceptional Enriques surfaces, as  developed by Ekedahl and Shepherd-Barron.}



\TitleInFrench{Le quotient fini unipotent maximal, torsion inhabituelle dans les solides de Fano et surfaces d'Enriques exceptionnelles.}

\AbstractInFrench{Nous introduisons et étudions le quotient fini unipotent maximal pour les schémas en groupes algébriques en caractéristique positive. Appliqué aux schémas de Picard, ce quotient encode la torsion inhabituelle. Nous construisons des solides de Fano intègres qui font effectivement apparaître cette torsion inhabituelle. L'existence de ces solides est surprenante, compte tenu du fait que la torsion inhabituelle s'annule pour les surfaces de Del Pezzo. Notre construction  repose sur la théorie des surfaces d'Enriques d\'evelopp\'ee par Ekedahl et Shepherd-Barron.}

\acknowledgement{
This research was  conducted in the framework of the   research training group
\emph{GRK 2240: Algebro-geometric Methods in Algebra, Arithmetic and Topology};
we wish to thank the Deutsche Forschungsgemeinschaft for financial support. 
The first-named author is currently funded by the \emph{Fondation Math\'ematique Jacques Hadamard.}}

\begin{document}


\removeabove{10pt}
\removebetween{10pt}
\removebelow{10pt}

\maketitle

\begin{prelims}

\DisplayAbstractInEnglish

\bigskip

\DisplayKeyWords

\medskip

\DisplayMSCclass

\bigskip

\languagesection{Fran\c{c}ais}

\bigskip

\DisplayTitleInFrench

\medskip

\DisplayAbstractInFrench

\end{prelims}


\newpage

\setcounter{tocdepth}{1} 

\tableofcontents


\section*{Introduction}
\mylabel{Introduction}

In algebraic geometry over ground fields $k$ of characteristic $p>0$,
unusual behavior of certain algebraic schemes   is often 
reflected by the structure of   unipotent torsion  originating from  the Picard group.
For example, an elliptic curve $E$ is \emph{supersingular} if and only if the kernel
$E[p]$  for multiplication-by-$p$ is unipotent.
An even more instructive case are Enriques surfaces $Y$, which  have $c_1=0$ and $b_2=10$.
Then the Picard scheme $P=\Pic^\tau_{Y/k}$ of numerically trivial invertible sheaves has
order two. In characteristic $p=2$, this gives the three possibilities.
In case $P=\mu_2$, the Enriques surface is called ordinary, and behaves like in characteristic zero.
Otherwise, we have $P=\ZZ/2\ZZ$ or $P=\alpha_2$, which is a unipotent group scheme, and $Y$ is a \emph{simply-connected} Enriques surfaces.
Their geometry and deformation theory is    more difficult to understand.
The crucial difference to the case of elliptic curves is that the unipotent torsion 
can be regarded as a \emph{quotient object}, and not only as a subobject, which makes it more ``unusual''.

The first  goal of this paper is to introduce a general measure  for unipotent torsion,
the \emph{maximal   finite unipotent quotient} $\Upsilon_{Y/k}=\Upsilon_P$ of the algebraic group scheme $P=\Pic^\tau_{Y/k}$. 
Over algebraically closed fields,
this comprises the $p$-primary torsion part of the N\'eron--Severi group $\NS(Y)$ 
and the  part of the local group scheme $P^0/P^0_\red$ whose Cartier dual is also local. Such finite group schemes
are often called \emph{local-local}.  
This actually extends to  algebraic group scheme that  are not necessarily commutative.
Our approach builds on the work of Brion \cite{Brion2017}. 
It turns out that  $\Upsilon_{Y/k}$ is useful in various situations. For example, it easily
explains that the reduced part  of an algebraic group scheme is not necessarily a subgroup scheme.
Of course, this   may only happen over imperfect fields.

The second goal is to construct  Fano varieties  whose Picard scheme actually contains such
unipotent torsion.
Roughly speaking, a \emph{Fano variety} is a   Gorenstein scheme $Y$
that is proper and equi-dimensional, and whose dualizing sheaf $\omega_Y$ is anti-ample. 
This notion indeed goes back to Fano \cite{Fano1931}. Usually, one also assumes that $Y$ is integral.
We write $n=\dim(Y)$ for the dimension. 
The structure and classification of Fano varieties is an interesting subject of its own. 
Fano varieties play an important role  in representation theory, because proper homogeneous spaces $Y=G/H$
for linear   groups schemes in characteristic $p=0$ are Fano varieties.
Moreover, they are crucial   for the minimal model program,
because they arise as generic fibers $Y=X_\eta$ in Mori fibrations $X\ra B$.

Suppose that $Y$ is a smooth Fano variety in characteristic zero.
Then  
$$
H^i(Y,\O_Y)= H^i(Y,\Omega^n_{Y/k}\otimes \omega_Y^{\otimes-1})=0
$$
for $i>0$, by Kodaira--Akizuki--Nakano Vanishing   (see \cite{Deligne-Illusie1987} for an algebraic proof).
In particular, the Lie algebra $H^1(Y,\O_Y)$ of the Picard scheme  $\Pic_{Y/k}$ vanishes,
so that the connected component $\Pic^0_{Y/k}$ is trivial. Moreover, the Euler characteristic is $\chi(\O_Y)=1$,
whence all connected \'etale coverings are trivial. Consequently, the algebraic fundamental
group   vanishes, and  the N\'eron--Severi group $\NS(Y)$
is torsion-free. Summing up, the group scheme  $\Pic^\tau_{Y/k}$ of numerically trivial invertible sheaves vanishes.
Its fundamental role was emphasized by Grothendieck \cite{FGAVI}.

It is a natural question to what extent the vanishing for   $\Pic^\tau_{Y/k}$
holds true in positive characteristics, or for singular Fano varieties.
Vanishing holds in dimension $n=1$ for integral Gorenstein curves with $\omega_Y^{\otimes -1}$ ample, basically by Riemann--Roch.
In dimension $n=2$, the classification of regular del Pezzo surfaces is independent of the characteristic,
and we again have vanishing. For normal del Pezzo surfaces, the second author observed
in \cite{Schroeer2001} that $h^1(\O_Y)=0$.
The situation becomes much more challenging for 
non-normal del Pezzo surfaces.
Here Reid \cite{Reid1994} constructed for each prime $p>0$ examples with $h^1(\O_Y)\neq 0$.
In \cite{Schroeer2007}, some normal locally factorial del Pezzo surfaces over imperfect fields 
in characteristic $p=2$ with $h^1(\O_Y)\neq 0$ were constructed, and even regular
examples exist, as shown by Maddock \cite{Maddock2016}. 
For smooth Fano threefolds, the vanishing of $\Pic^\tau_{Y/k}$ was  shown by Shepherd-Barron \cite{Shepherd-Barron1997}.
However, Cascini and Tanaka \cite{Cascini-Tanaka2016} constructed a klt Fano threefold in characteristic $p=2$ with $h^2(\O_Y) \neq 0$;
other examples with $p\ge 3$ were found by Bernasconi \cite{Bernasconi2019}.

Moreover, Tanaka \cite{Tanaka2016} constructed Mori fibrations  $X\ra B$ on threefolds in characteristic $p=2,3$
where the generic fiber $Y=X_\eta$ is a normal del Pezzo surface whose Picard group  contains elements of order $p=2,3$.
Very recently, Bernasconi and Tanaka \cite{Bernasconi-Tanaka2019} provided effective bounds for the torsion on del Pezzo--type 
surfaces over imperfect fields.

However, in all the above
examples the maximal unipotent finite quotient $\Upsilon_{Y/k}$ of $\Pic^\tau_{Y/k}$ still
vanishes. We therefore regard $\Upsilon_{Y/k}$ as a measure for \emph{unusual torsion 
in Fano varieties}.
The main results of this paper are as follows.
First, we give general criteria for the vanishing of $\Upsilon_{Y/k}$,
which   relies on the theory of Bockstein operators:

\begin{maintheorem}
{\rm (see Thm.\ \ref{hasse--witt})}
We have $\Upsilon^0_{Y/k}=0$ provided the Frobenius map on $H^2(Y,\O_Y)$ has maximal Hasse--Witt rank.
\end{maintheorem}

From this  we deduce in Theorem \ref{vanishing surfaces} that $\Upsilon_{Y/k}=0$
for any reduced surface whose dualizing sheaf $\omega_Y$ is negative, in a suitable sense.
Note that this requires no further restrictions on the singularities.
It applies in particular to   del Pezzo surfaces.
In light of this, it is surprising that unusual torsion does appear in dimension $n\geq 3$:

\begin{maintheorem}
\text{\rm (see Thm.\ \ref{exotic threefold})}
There are integral Fano threefolds $Y$ in characteristic two such that 
$\Upsilon_{Y/k}=\Pic^\tau_{Y/k}$ is isomorphic to the group scheme  $\ZZ/2\ZZ$ or $\alpha_2$.
\end{maintheorem}

The construction relies on the theory of \emph{exceptional Enriques surfaces},
as developed by Ekedahl and Shepherd-Barron \cite{Ekedahl-Shepherd-Barron2004}.
These are simply-connected Enriques surfaces containing very strange configurations
of $(-2)$-curves. The geometry of such surfaces  was already studied in the monograph of Cossec and Dolgachev
(\cite{Cossec-Dolgachev1989}, Chapter III, \S3--4), although their existence was
established only later by Salomonsson \cite{Salomonsson2003}. 

The idea for our construction is rather simple: start with
an     exceptional Enriques surface $S$ that contains exactly ten $(-2)$-curves
with a certain dual graph. We then contract all but one of these curves,
creating a normal Enriques surface $Z$. Its normalized K3-like covering $Z'\ra Z$ turns out to be  a normal
del Pezzo surface with a unique singularity, which is a rational double point
of type $D_5$.  This relies on the classification of
normal del Pezzo surfaces, as explained  by Dolgachev \cite{Dolgachev2012}.
We   take the $\PP^1$-bundle $X=\PP(\shE)$ for the
locally free sheaf $\shE=\O_{Z'}\oplus\omega_{Z'}^{\otimes-1}$,
contract the negative section $E=\PP(\omega_{Z'}^{\otimes-1})$, and 
denormalize along the positive section $Z'=\PP(\O_{Z})$
with respect to the purely inseparable double covering $\nu:Z'\ra Z$.
The resulting scheme $Y=X\amalg_Z Z'$ is an integral Fano threefold
with Euler characteristic,  degree   and  index  
$$
\chi(\O_Y)=1,\quad \deg(Y)=4\quadand \ind(Y)=1.
$$ 
A detailed analysis of the normalized K3-like covering $\nu:Z'\ra Z$ reveals that
the morphism is flat, except over a rational double point  $a\in Z$ of type $E_8$.
From this we infer that our integral Fano threefold $Y$ has invertible dualizing sheaf
and satisfies Serre  Condition $(S_2)$, but \emph{fails to be Cohen--Macaulay}. 
Such schemes might be called \emph{quasi-Fano varieties}, and deserve further study.

Note that recently Totaro \cite{Totaro2019} constructed three-dimensional terminal singularities in positive characteristics that
are not Cohen--Macaulay, and further examples were constructed by Yasuda \cite{Yasuda2019}.
In turn, the above examples of non-normal Fano threefolds $Y$ may admit twisted forms $Y'$   whose
local rings are normal $\QQ$-factorial klt singularities, and thus could occur as generic fibers 
in Mori fiber   spaces. See \cite{Fanelli-Schroeer2020} for our analysis
of non-normal del Pezzo surfaces having   twisted forms whose local rings are regular.

\medskip
The paper is organized as follows:
In Section \ref{Maximal unipotent} we discuss  various smallest normal subgroup schemes
$N\subset G$ whose quotients   have certain properties,
introduce  the maximal finite unipotent quotient $\Upsilon_G=G/N$, and establish
its basic properties.
We apply this in Section \ref{Picard scheme} to  Picard schemes $G=\Pic^\tau_{Y/k}$,
and give general criteria for the vanishing of $\Upsilon_{Y/k}=\Upsilon_G$.
In Section \ref{Surfaces} we show that $\Upsilon_{Y/k}$ vanishes for
reduced surfaces with negative  dualizing sheaf. 
In Section \ref{Enriques del Pezzo}  we recall various notions of Enriques surfaces,
and describe  how   del Pezzo surfaces arise as normalization of K3-like
coverings.
The following  Section \ref{Exceptional enriques} contains a detailed analysis for
the  case of  exceptional Enriques surfaces of type $T_{2,3,7}$.
In Section \ref{Cones} we discuss the cone construction and its Gorenstein properties,
and show how it leads to new Fano varieties.
The final Section \ref{Fano threefolds} contains the construction 
of integral Fano threefolds whose maximal finite unipotent quotient $\Upsilon_{Y/k}$
is non-trivial.
\vskip\baselineskip
\noindent\textbf{Acknowledgements:} We like to thank Fabio Bernasconi, Michel Brion, Igor Dolgachev, Hiromu Tanaka and the referee for many valuable suggestions.

\section{The maximal finite unipotent quotient}
\mylabel{Maximal unipotent}

In order to study unusual torsion in Picard groups, we shall  introduce
the maximal finite unipotent quotient for general algebraic group schemes, which are not
necessarily commutative. These results seem to be of independent interest, and
mainly rely on the theory of algebraic groups.

Let $k$ be a ground field of characteristic $p\geq 0$. 
An \emph{algebraic group scheme} is a group scheme $G$ where the structure morphism $G\ra\Spec(k)$ is  of finite type.
Note that the underlying scheme is automatically separated.
We say that $G$ is \emph{finite} if the structure morphism is finite. Then the \emph{order}  is   defined as 
$\ord(G)=h^0(\O_G)=\dim_k H^0(G,\O_G)$.
One says that  $G$ is   \emph{of multiplicative type} if  the base-change $G\otimes k^\alg$
to the algebraic closure is isomorphic to $\Spec k^\alg[M]$, where $M$ is a finitely generated abelian group.
In other words, $G$ is commutative and  the base-change has a filtration whose subquotients are isomorphic to the multiplicative
group $\GG_m$, or the constant groups $\mu_l=\ZZ/l\ZZ$ for some prime $l\neq p$, or the
local group scheme $\mu_p$. One also says that $G$ is \emph{multiplicative}.
The algebraic group  $G$ is  called \emph{unipotent}
if  $G\otimes k^\alg$   admits a filtration whose subquotients 
are isomorphic to subgroup schemes of  the additive group $\GG_a$.
In characteristic $p>0$, this means that after refinement the subquotients are $\GG_a$, 
the constant group $\ZZ/p\ZZ$ or the local group scheme $\alpha_p$.
Note that unipotent group schemes are not necessarily commutative.
For more on these   notions, see \cite[expos\'e IX and expos\'e XVII]{SGA3b}.
 
\begin{theorem}
\mylabel{maximal quotients}
For each of the following properties {\rm(i)-(vi)} and each  algebraic group scheme $G$,
there is a smallest normal subgroup scheme $N\subset G$ such that the quotient $G/N$ has the
property in question:
\begin{multicols}{2}
\begin{enumerate}
\item \'etale;
\item affine;
\item proper;
\item finite;
\item finite and multiplicative;
\item finite and unipotent. 
\end{enumerate}
\end{multicols}
\noindent 
Moreover, the   subgroup schemes $N\subset G$  and  the quotients $G/N$ commute with separable field extensions $k\subset k'$.
\end{theorem}

\proof
By \cite[Expos\'e $\rm{IV_A}$, Section 2]{SGA3a}, the connected component $G^0\subset G$ of the origin is 
a normal subgroup scheme such that $G/G^0$ is \'etale, and it is indeed the smallest.
According to Brion's analysis (\cite[Theorem  1 and Theorem 2]{Brion2017}), there are the smallest normal subgroup schemes $N_1,N_2\subset G$ such that
the resulting quotients are affine and proper, respectively.  Let $N\subset G$ be the normal subgroup scheme
generated by $N_1$ and $N_2$. Then the resulting quotient $G/N$ is both proper and affine, hence finite.
Moreover, every homomorphism $G\ra K$ into some finite group scheme   contains $N_1$ and $N_2$ in
its kernel. It follows that $N$ is the desired smallest normal subgroup scheme with finite quotient.
This settles the cases (i)--(iv).

Consider the ordered family of normal subgroup schemes 
$H_\lambda\subset G$, $\lambda\in L$ whose quotients $G/H_\lambda$
are finite and unipotent. The group scheme $G$ itself belongs to this family, and each member
contains $N$.
We first check that for any two  members $H_\lambda$ and $H_\mu$, the intersection 
$K=H_\lambda\cap H_\mu$ also belongs to the family.
We have an exact sequence
$$
0\lra H_\lambda/K \lra G/K\lra G/H_\lambda\lra 0.
$$
By the Isomorphism Theorem, the term on the left is isomorphic to $(H_\lambda\cdot H_\mu)/H_\mu$, which is contained in $G/H_\mu$.
The latter is finite and unipotent, so the same holds for the subgroup scheme $(H_\lambda\cdot H_\mu)/H_\mu$ and
the extension $G/K$.

Seeking a contradiction, we assume that there is no smallest member.
Since the family is filtered, this means that it contains an infinite    descending sequence
$H_0\supsetneqq H_1\supsetneqq\ldots$  such that  the quotients   $U_n=G/H_n$ have unbounded orders.
On the other hand, all of them are quotients of the finite group scheme $G/N$, hence the orders are bounded,
contradiction. Hence, there is a smallest normal subgroup scheme whose quotient is finite and unipotent.
This settles (vi).
The argument for   (v) is similar and left to the reader.

We now prove the second part of the assertion: 
let $k\subset k'$ be a separable extension,     
and   $N'\subset G'\otimes k'$ be the smallest subgroup scheme over $k'$ 
such that the quotient has  the property $\mathcal{P}$ in question. This gives an inclusion $N'\subset N\otimes k'$,
and we have to verify that it is an equality.
Suppose first that $k\subset k'$ is algebraic. By fpqc descent, it suffices to check $N'=N\otimes k'$
after enlarging the field extension. It thus suffices to treat the case that  $k\subset k'$ is Galois, 
with Galois group $\Gamma=\Gal(k'/k)$.  Then for each element $\sigma\in \Gamma$ we have
$\sigma(N')=N'$, by the uniqueness of $N'\subset G\otimes k'$. Galois descent gives  a closed
subscheme $N_0\subset N$ with $N_0\otimes k'=N'$.
This subscheme is a subgroup scheme and normal in $G$, and the base-change $G/N_0\otimes k'=(G\otimes k')/N'$
has property $\mathcal{P}$. For each of the cases (i)--(vi),
this implies that  $G/N_0$ has property $\mathcal{P}$. This gives $N_0=N$, and in turn 
the desired equality $N'=N\otimes k'$.
 
Now write $k'=\bigcup k_\lambda$ as the filtered union of finitely generated subextensions.
The structure morphism $G'\ra\Spec(k')$ is of finite presentation. Hence, \cite[Theorem 8.8.2]{EGAIVc} ensures that
there is some index $\lambda$ so that the closed subscheme $N'$ is the base-change of some closed subscheme 
$N_\lambda\subset G\otimes k_\lambda$.
Moreover,    this subscheme is a normal subgroup scheme, and we have $N_\lambda\subset N\otimes k_\lambda$.
This reduces our problem to the case that $k'$ is finitely generated. 
Choose an integral affine scheme $S$ of finite type 
with function field $\kappa(\eta)=k'$.
Since $k\subset k'$ is separable, the scheme $S$ is geometrically reduced.
Passing to some dense open  set,  we may assume that $S$ is smooth.

Consider the relative group scheme $G_S$, with generic fiber $G_\eta=G\otimes k'$.
Seeking a contradiction, we assume that $N'\subset N\otimes k'$ is not an equality.
Then there is a homomorphism $f_\eta\colon G_\eta\ra H_\eta$ to some algebraic group scheme $H_\eta$ 
having   property $\mathcal{P}$, such that $N_\eta$ is not contained in the kernel.
Again by \cite[th\'eor\`eme~8.8.2]{EGAIVc}, there is a dense open  set $U\subset S$,
a relative group scheme $H_U$ where the structure morphism $H_U\ra U$ is   of finite type,
and a homomorphism $f_U \colon G_U\ra H_U$ inducing $f_\eta$. 
If the scheme $H_\eta$ is affine, proper or finite, we may assume that the respective property holds
for the morphism $H_U\ra U$, by \cite[th\'eor\`eme~8.10.5]{EGAIVc}.
Moreover, if $H_\eta$ is finite and unipotent, we may assume that there is a finite \'etale
covering $U'\ra U$ so that $H\times_UU'$ is a successive extension of $\ZZ/p\ZZ$ and $\alpha_p$,
such that all fibers $H_U\ra U$ are finite and unipotent.
The situation for properties (i) and (v) is similar.
Summing up, we may assume that all fibers of $H_U\ra U$ have the property $\mathcal{P}$ in question.
Passing to a dense open   set again, we reduce to  $U=S$.
Consider  the $S$-scheme
$$
K=\Kernel(f|N)=N\times_{G} \{e_S\}.
$$
This  is a relative group scheme, and the 
structure morphism $g \colon K\ra S$ is of finite type.
The generic fiber  is equi-dimensional, say of dimension $n=\dim(K_\eta)$.
By \cite[expos\'e ${\rm VI_B}$, proposition 4.1]{SGA3a},
the set of   points $a\in S$ with  $\dim(K_a)=n$  is constructible.
Replacing $S$ by some dense open set, we may assume that all fibers $K_a$ are $n$-dimensional.
By Bertini's Theorem (\cite[th\'eor\`eme~6.3]{Jouanolou1983}),  there are closed points $a\in S$ such that
the finite extension $k\subset\kappa(a)$ is separable. Since the fiber $H_a$ has property $\mathcal{P}$ in question,
the kernel $K_a$ is trivial, and thus $n=0$. In turn, the morphism $g\colon K\ra S$ is quasi-finite.
By Zariski's Main Theorem, there is a closed embedding $K\subset X$ into some finite 
$S$-scheme $X$, with $K_\eta=X_\eta$. After replacing $S$ by a dense open  set, we may assume
that $g\colon K\ra S$ is finite, and  furthermore  flat, say of degree $d=\deg(K/S)$. 
Looking again at   fibers $K_a$, we see that $d=1$. In turn, $K_\eta$ is trivial,
contradiction. 
\qed

\medskip
Let us say  that $G/N$ is  \emph{the maximal quotient}  with the   property  $\mathcal{P}$ in question.
By construction, any homomorphism $G\ra H$ into some algebraic group scheme with property $\mathcal{P}$ uniquely factors over $G/N$.
Thus $G/N$ is functorial in $G$, and the    functor $G\mapsto G/N$ is the left adjoint for the inclusion 
$H\mapsto H$ of the category of  algebraic group schemes with property $\mathcal{P}$
into the category of all algebraic group schemes.

Given a field extension  $k\subset k'$, we set $G'=G\otimes k'$ and let $N'\subset G'$ be the
normal subgroup scheme giving the maximal quotient over $k'$. 
This gives an inclusion   $N'\subset N\otimes k'$, and a resulting \emph{base-change map}
$G'/N'\ra (G/N)\otimes k' $. It may or may not be an isomorphism, as we shall see below.
However, we have the following immediate fact:

\begin{proposition}
\mylabel{base-change vanish}
In the above situation, the base-change map is an epimorphism. In particular, the algebraic group scheme
$G/N$ vanishes if $G'/N'$  vanishes.
\end{proposition}

The maximal affine quotient is indeed the \emph{affine hull} in the sense of scheme theory, and written as
$G^\aff=\Spec \Gamma(G,\O_G)$. The kernel $N$ for the homomorphism $G\ra G^\aff$  is \emph{anti-affine},
which means that the inclusion $k\subset H^0(N,\O_N)$ is an equality (\cite[Chapter III, \S3, Theorem 8.2]{Demazure-Gabriel1970}). 
This notion was introduced by Brion \cite{Brion2009}, and implies 
 that $N$ is \emph{semi-abelian}, and in particular smooth and connected (\cite[Proposition 5.5.1]{Brion2017}).
The maximal \'etale quotient is usually denoted by $\Phi_G=G/G^0$, and called the \emph{group scheme
of   components}.  
Both $G^\aff$ and $\Phi_G$ actually commute with arbitrary field extensions $k\subset k'$.

Throughout this paper  we are mainly interested in the \emph{maximal finite unipotent quotient}
$$
\Upsilon_G=G/N.
$$
Any subgroup scheme of $G$ that is multiplicative, or smooth and connected,   or   merely integral vanishes in $\Upsilon_G$.
In fact, the following holds:

\begin{lemma}
\mylabel{zero to unipotent}
Suppose the algebraic group scheme $G$ has a filtration  whose subquotients are generated by integral group schemes
and group schemes of multiplicative type. Then all homomorphisms $f\colon G\ra U$ into finite unipotent group schemes $U$ are trivial.
\end{lemma}

\newcommand{\Image}{\operatorname{Im}}
\proof
By induction on the length of the filtration, it suffices to treat the two cases that $G$ is integral,
or of multiplicative type. In the latter case, the statement follows from \cite[expos\'e~XVII, proposition~2.4]{SGA3b}.
Now suppose that $G$ is integral. Replacing $U$ by the image $\Image(f)=G/\Kernel(f)$,
we may assume that  $f\colon G\ra U$ is surjective and schematically dominant. Setting $A=H^0(G,\O_G)$ and $R=H^0(U,\O_U)$,
we see that the canonical map $R\ra A$ is injective. These rings are integral, because the scheme $G$ is integral.
Moreover,  $R$ is an Artin ring, because the group scheme  $U$ is finite. The neutral element $e\in U$
shows that the residue field is $R/\maxid_R=k$.
In turn, we have  $R=k$, thus $U=\Spec(k)$ is trivial.
\qed

\medskip
Each algebraic group scheme $G$ yields a \emph{Galois representation} $G(k^\sep)$ of the Galois group $\Gamma=\Gal(k^\sep/k)$
on the abstract  group   $G(k^\sep)$, which might be infinitely generated. 
However, this construction yields an equivalence between the category
of algebraic group schemes that are \'etale and  continuous Galois representations on finite groups
(\cite[expos\'e~VIII, proposition~2.1]{SGA4b}).
Note also that for finite group schemes $H$, the group scheme $H^0$ is local,
and the   finite universal homeomorphism 
$H\ra\Phi_H$ yields an equality $H(k')=\Phi_H(k')$ for all field extensions $k\subset k'$.
This applies in particular to the maximal finite unipotent quotient $H=\Upsilon_G$.
The following observation thus computes its  group scheme of components:

\begin{proposition}
\mylabel{galois module}
The Galois module   $\Upsilon_G(k^\sep)$ is the quotient of the finite group  $\Phi_G(k^\sep)$
by the subgroup generated by all $l$-Sylow groups for the primes $l\neq p$.
\end{proposition}
 
\proof
Base-changing to $k^\sep$, we are reduced to the case that the ground field $k$ is separably closed.
Write $\Phi'_G$ for the quotient of the constant group scheme  $\Phi_G$ by the subgroup generated by the Sylow-$l$-subgroups,
and set $U=\Upsilon_G$. Then $\Phi_U(k)=U(k)$ is a finite $p$-group.
It follows that the canonical surjection $\Phi_G\ra\Phi_U$ factors over $\Phi'_G$.
The resulting map $f\colon \Phi'_G\ra\Phi_U$ admits a section, by the universal properties
of   $\Upsilon_G$ and $\Phi_U$. This section $s\colon \Phi_U\ra\Phi'_G$ is surjective, because $G\ra\Phi'_G$ is surjective.
In turn, $f$ and $s$ are   inverse to each other, and give the desired identification $\Phi'_G(k)=\Phi_U(k)=\Upsilon_G(k)$.
\qed

\medskip
It is much more difficult to understand the component of the origin $\Upsilon^0_G$,
because the base-change map $\Upsilon_{G\otimes k'}\ra \Upsilon_G\otimes k'$ might fail to be 
an isomorphism, for inseparable extensions $k'$. Here is a typical example, over imperfect fields $k$:

\begin{proposition}
\mylabel{not base-change}
Suppose  $G$ is an integral group scheme such that some base-change
$G'=G\otimes k'$  becomes isomorphic to $\GG_a\oplus\alpha_p$.
Then we have $\Upsilon_G=0$, whereas  $\Upsilon_{G'}=\alpha_p$. 
\end{proposition}

\proof
According to  Proposition \ref{zero to unipotent}, the quotient  $\Upsilon_G$ is trivial,
whereas $G'\ra\Upsilon_{G'}$ is the projection onto the second factor.
\qed

\medskip
For example, $G$ could be the kernel for the homomorphism $h\colon \GG_a\times\GG_a\ra \GG_a$
given by the additive map $(x,y)\longmapsto x^p+ty^p$, for some scalar $t\in k$ that is not a $p$-th power.
Note that  certain fiber products  $G=\GG_a\times_{\GG_a}\GG_a$ were systematically studied by Russel \cite{Russel1970}
to describe twisted forms of the additive group.
In our concrete example the  compactification in $\PP^2$ is the Fermat curve with homogeneous equation $X^p+tY^p+Z^p=0$,
as analyzed in \cite[\S4]{Schroeer2010}. The scheme $G$ is integral, with singular locus $G(k)=\{0\}$,
and its normalization is the affine line over the height-one extension $E=k(t^{1/p})$.
It has another peculiar feature:

\begin{proposition}
\mylabel{not subgroup}
Notation as above.
For the connected scheme $H=G\times G$, the closed subscheme $H_\red\subset H$ is not a subgroup scheme.
\end{proposition}

\proof
This already appeared in \cite[expos\'e~VI, exemple~1.3.2]{SGA3a}. Let us give an independent short  proof
based on our maximal finite unipotent quotient: Suppose $H_\red $ were a subgroup scheme. Then $H/H_\red$ is a local
unipotent group scheme, so the projection factors over $\Upsilon_H=\Upsilon_G\times\Upsilon_G=0$,
consequently $H_\red=H$. On the other hand, the normalization $\AA^1_E\ra G$ shows that $H$ is birational
to the affine plane over the non-reduced ring $E\otimes E=E[t]/(t^p)$, hence $H_\red\subsetneqq H$,
contradiction.
\qed

\medskip
If the reduced part $G_\red\subset G$ of an algebraic group scheme   is \emph{geometrically reduced},
then the product  $G_\red\times G_\red$ remains  reduced, and the group multiplication factors
over $G_\red$. In turn, the closed subscheme $G_\red\subset G$ is a subgroup scheme.
This frequently fails, as we saw above. 
Other examples for this behavior are  the  non-split extensions $0\ra \alpha_p\ra G\ra \ZZ/p\ZZ\ra 0$,
where the fiber $\pr^{-1}(a)$ for the projection $\pr\colon G\ra\ZZ/p\ZZ$ is reduced if and only if $a\neq 0$.
Such extensions exist over imperfect fields: the  abelian group of all \emph{central extensions} contains the flat cohomology group
$H^1(\Spec(k),\alpha_p)=k/k^p$ as a subgroup, by \cite[Chapter III, \S6, Proposition in 3.5]{Demazure-Gabriel1970}.
Also note that   if $G_\red$ is a subgroup  scheme, it need not be normal,  for example in  semidirect products like $G=\alpha_p\rtimes\GG_m$.

The situation simplifies somewhat  for commutative group schemes.
Recall that if $G$ is commutative and affine, then there is  a \emph{maximal multiplicative subgroup scheme}
$G^{\mult}\subset G$, and the quotient is unipotent
(\cite[Chapter IV, \S3, Theorem 1.1]{Demazure-Gabriel1970}). If the resulting quotient is finite,  
Lemma \ref{zero to unipotent}  gives  $\Upsilon_G=G/G^{\mult}$.
In general, we get:

\begin{proposition}
\mylabel{component orgin}
Suppose that $G$ is commutative, that the  scheme $G^0_\red$ is geometrically reduced, and
that the projection $G\ra\Phi_G$ admits a section. For the local group scheme $L=G^0/G^0_\red$,
we get an identification $\Upsilon^0_G=L/L^{\mult}$.
\end{proposition}

\proof
The section gives a decomposition of commutative algebraic groups $G=G'\oplus G''$, where the first factor is
connected and the second factor is \'etale. In turn, we have $\Upsilon_G=\Upsilon_{G'}\times\Upsilon_{G''}$.
From the universal property we infer that $\Upsilon_{G''}$ is \'etale, so we may assume from the start that $G$ is connected,
and have to show that the projection $G\ra G/G_\red=L$ induces an isomorphism $\Upsilon_G\ra\Upsilon_L=L/L^{\mult}$.
Since $G_\red$ is geometrically reduced, the inclusion $G_\red\subset G$ is a subgroup scheme,
which is smooth and connected. It lies in the kernel of any homomorphism $G\ra U$ to some finite unipotent
scheme $U$, by Lemma \ref{zero to unipotent}. From the universal properties we infer that  $\Upsilon_G\ra\Upsilon_L$ is an isomorphism.
\qed
 
\medskip
This leads to the following structure result:

\begin{theorem}
\mylabel{structure upsilon}
If $G$ is commutative and $k$ is perfect, we get an identification
$$\Upsilon_G=L/L^{\mult}\times \Phi_G[p^\infty],$$ with the local group scheme  $L=G^0/G^0_\red$ 
and the $p$-primary part $\Phi_G[p^\infty]\subset\Phi_G$.
Moreover, the kernel $N$ for the projection  $G\ra\Upsilon_G$ has a   three-step filtration 
with $N/N_2$     multiplicative, $N_2/N_1$   smooth unipotent, and $N_1$ anti-affine.
\end{theorem}

\proof
Since $\Phi_G$ is commutative, the quotient of $\Phi_G(k^\sep)$ by the subgroup generated by the Sylow-$l$-groups
gets identified with the $p$-primary torsion in $\Phi_G(k^\sep)$.
From Proposition \ref{galois module} we infer that $\Phi_G[p^\infty]$ is the group scheme of components for $\Upsilon_G$. 
Since $k$ is perfect, the scheme $G^0_\red$ is geometrically reduced, and Proposition
\ref{component orgin} gives $\Upsilon^0_G=L/L^\mult$.   
Moreover, the reduced part of $\Upsilon_G$ yields a section for 
the group scheme of components, and we get the decomposition $\Upsilon_G=L/L^{\mult}\times \Phi_G[p^\infty]$.
Applying the Five Lemma to the commutative diagram
$$
\begin{CD}
0	@>>>	G^0	@>>>	G	@>>>	\Phi_G	@>>> 0\\
@.		@VVV		@VVV		@VVV\\
0	@>>>	L/L^\mult @>>>	\Upsilon_G @>>>	\Phi_G[p^\infty]	@>>>0,
\end{CD}
$$
we see that the kernel $N$ for $G\ra \Upsilon_G$ has a  filtration $F_i\subset N$ with
$$
F_5=N\quadand F_4=\Kernel(G^0\ra L/L^\mult)\quadand F_3=\Kernel(G^0\ra L)=G^0_\red.
$$
Then $F_5/F_4$ is the sum of the $l$-primary parts in $\Phi_G$ for the primes $l\neq p$,
and $F_4/F_3=L^\mult$.
The kernel $F_1$ for the affinization $F_3\ra F_3^\aff$ is anti-affine.
Moreover, the  multiplicative part of $F_3^\aff$ has a unique complement, because $k$ is perfect,
and this complement is smooth unipotent. Let $F_2\subset F_3$ be its preimage.

Summing up, we have constructed a five-step filtration $0=F_0\subset \ldots\subset F_5=N$.
Setting $N_2=F_2$ and $N_1=F_1$ we see that $N/N_1$ is multiplicative, $N_2/N_1$ is smooth unipotent,
and $N_1$ is anti-affine.
\qed

\medskip
Note that the  epimorphism $G\ra\Upsilon_G$ does not admit a section in general:
suppose that $N$ is either the additive group, or a supersingular elliptic curve.
With respect to the scalar multiplication of $\End(\alpha_p)=k$,
the abelian group $\Ext^1(\alpha_p,N)$ becomes  a one-dimensional vector space,
provided that $k$ is algebraically closed
(\cite[table on page II.14-2]{Oort1966}).
In the ensuing non-split extensions $0\ra N\ra G\ra \alpha_p\ra 0$, the projection
coincides with $G\ra\Upsilon_G$, according to Lemma \ref{zero to unipotent}.

\section{Picard scheme and Bockstein operators}
\mylabel{Picard scheme}

Let $k$ be a ground field of characteristic $p>0$, and $Y$ be a proper scheme.
Then the \emph{Picard scheme} $P=\Pic_{Y/k}$ exists, and this group scheme is locally of finite type
(\cite[expos\'e~XII, corollaire~1.5 and expos\'e~XIII, proposition~3.2]{SGA6}).
The Galois module for the \'etale group scheme $\Phi_P=P/P^0$  is the \emph{N\'eron--Severi group}
$\Phi_P(k^\sep)=\NS(Y\otimes_kk^\alg)$, which is finitely generated. 
In particular, the torsion part in $\Phi_P$ is a finite group scheme, hence    its inverse image $G=\Pic^\tau_{Y/k}$ 
in the Picard scheme is an  algebraic group scheme. We now  consider the maximal finite unipotent quotient
$$
\Upsilon_{Y/k} =  \Upsilon_G=\Upsilon_{\Pic^\tau_{Y/k}},
$$
and  regard this in many situations 
as a measure for   \emph{unusual} behavior of torsion in the  Picard scheme.  From Proposition \ref{galois module} we get:

\begin{proposition}
\mylabel{galois module upsilon}
The Galois module $\Upsilon_{Y/k}(k^\sep)$ is the $p$-torsion part of the
N\'eron--Severi group $\NS(Y\otimes k^\alg)$.  
\end{proposition}

\newcommand{\HW}{{\text{\rm HW}}}
We now seek to   understand the   component of the origin $\Upsilon_{Y/k}^0$.
For this the Frobenius map  $f\mapsto f^p$ on the structure sheaf $\O_Y$ is crucial. 
This map is additive, and becomes $k$-linear when one re-defines scalar multiplication 
on the range as $\lambda \cdot f=\lambda^pf$. Such additive maps are called \emph{$p$-linear}.
We now consider the induced $p$-linear maps on  the cohomology groups  $H^i(Y,\O_Y)$.

To understand this   better, suppose we have an arbitrary finite-dimensional $k$-vector space
$V$, together with a $p$-linear map $f\colon V\ra V$. 
Choose a  basis $a_1,\ldots,a_n\in V$. Then $f$ is determined by the images $f(a_j)$,
and the expansion $f(a_j) = \sum \lambda_{ij}a_i$ gives a matrix $A=(\lambda_{ij})\in\Mat_n(k)$.
A different basis $b_1,\ldots,b_n\in V$  yields another matrix  $B=SAT$, where 
the base-change matrix $S=(\sigma_{ij})$ is defined by $a_j=\sum\sigma_{ij}b_i$, and $T$ is obtained
from the inverse of $S$ by applying Frobenius to the entries.
In turn, the rank of    $A$ depends only on $f$. Let us  call this integer the
\emph{Hasse--Witt rank} $\rank_\HW(f)\geq 0$ of the $p$-linear map $f\colon V\ra V$.
Note also that one may define the \emph{Hasse--Witt determinant}
$\det_\HW(f)=\det(A)$ as a   class in the monoid $k/k^{\times (p-1)}$.
Furthermore, we may regard the datum $(V,f)$ as a \emph{left  module} over the associative ring $k[F]$, where the 
relations $F\lambda = \lambda^pF$ hold, by setting $F\cdot a=f(a)$.
The  Hasse--Witt rank and determinant then become    invariants of this  module.
All these considerations go  back to Hasse and Witt \cite{Hasse-Witt1936},
who studied  $V=H^1(C,\O_C)$ for a smooth algebraic curve $C$ over an algebraically closed field $k$.
Compare also the recent discussion of Achter and Howe \cite{Achter-Howe2019} for a discussion of 
historical developments, and 
widespread inaccuracies in the literature.

Given a field extension $k\subset k'$, we see that there is a  unique $p$-linear extension $f'$ of $f$
to $V'=V\otimes k$, given by $f(a_j\otimes \lambda)=\lambda^pf(a_j)$.
Obviously, the $p$-linear maps $f$ and $f'$ have the same Hasse--Witt rank.
If $k$ is perfect, the subgroup $U=f(V)$ inside $V$ is actually  a vector subspace with respect
to the original scalar multiplication, and we have $\rank_\HW(f)=\dim_k(U)$.
If $\rank_\HW(f)=\dim(V)$, we say that  $f$ has \emph{maximal Hasse--Witt rank}.
This means that $f(V)$ generates the vector space $V$ with respect to the original
scalar multiplication. Equivalently,  for   some and hence all perfect field extensions $k\subset k'$
the $p$-linear extension $f'\colon V'\ra V'$ is bijective.
We now come to the main result of this section:
 
\begin{theorem}
\mylabel{hasse--witt}
We have $\Upsilon_{Y/k}^0=0$ provided 
the Frobenius map on the second cohomology group $H^2(Y,\O_Y)$ has maximal Hasse--Witt rank.
\end{theorem}

The proof is deferred to the end of this section. It 
relies on   Bockstein operators, a theory  introduced by Serre \cite{Serre1958}, which we like
to discuss first.
Write $W_m(k)$ be the ring of  \emph{Witt vectors}
$(\lambda_0,\ldots,\lambda_{m-1})$ of length $m$.
This is a ring endowed with two commuting additive maps  \emph{Frobenius} $F$ and   \emph{Verschiebung} (``shift'') $V$,
given by the formula
$$
F(\lambda_0,\ldots,\lambda_{m-1}) = (\lambda_0^p,\ldots,\lambda_{m-1}^p)\quadand
V(\lambda_0,\ldots,\lambda_{m-1}) = (0,\lambda_0,\ldots,\lambda_{m-2}).
$$
We refer to \cite[chapitre~9, \S1]{AC8-9} for the general theory of  Witt vectors.
The canonical projection $W_m(k)\ra W_{m-n}(k)$ is a homomorphism of rings, whose kernel we denote by
$$
V_m^n(k)=V^{m-n}W_m(k).
$$
Note that this kernel has length $m-n$, and  is stable under Frobenius, by the relation $FV=VF$.
Likewise, we have a short exact sequence of abelian sheaves
$$
0\lra V^n_m(\O_Y)\lra W_m(\O_Y)\lra  W_{m-n}(\O_Y)\lra 0,
$$
where the maps are $W_m(k)$-linear and compatible  with Frobenius. To simplify notation, 
write $W_m=W_m(\O_Y)$ and $V_m^n=V_m^n(\O_Y)$,
such that $W_1=\O_Y$.  
Combining the   long  exact sequences for the short exact sequences
$0\ra V^r_{r+1}\ra W_{r+1}\ra W_r\ra 0$ and $ 0\ra V^1_r\ra W_r\ra \O_Y\ra 0$, 
we get for each $i\geq 0$ a  commutative diagram
$$
\begin{tikzcd}[]
				& H^{i}(W_{r+1})\ar[d]\ar[dr]\\
H^{i}(V^1_r)\ar[r]\ar[dr]	& H^{i}(W_r)\ar[r]\ar[d]		& H^{i}(\O_Y)\\
				& H^{i+1}(V^r_{r+1})
\end{tikzcd}
$$
with exact row and column. This gives a canonical $W_r(k)$-linear map
$$
\Image(H^{i}(W_r)\ra H^{i}(\O_Y))\;\stackrel{\beta_r}{\lra} \;\Cokernel(H^{i}(V^1_r)\ra H^{i+1}(V^r_{r+1}))
$$
called \emph{Bockstein operator}, by sending the image of   $x\in H^{i}(W_r)$ in $H^{i}(\O_Y)$
to  the image of $x$
in $H^{i+1}(V^r_{r+1})$ modulo the image of $H^{i}(V^1_r)$. Obviously, the kernel of $\beta_r$ comprises those cohomology classes in $H^{i}(\O_Y)$
that extend to $H^{i}(W_{r+1})$. These   form a decreasing sequence.
We write $H^i(Y,\O_Y)[\beta]$ for their common intersection, and call it the \emph{Bockstein kernel}.
By construction, this is a vector subspace of $H^i(Y,\O_Y)$ invariant under Frobenius. 

Now consider the case $i=1$. Then $\lieg=H^1(Y,\O_Y)$ is the \emph{Lie algebra} for the
Picard scheme, hence also for the algebraic group scheme $G=\Pic^0_{Y/k}$.
As such, it has an additional structure, namely the \emph{$p$-power map} $x\mapsto x^{[p]}$ obtained
from the $p$-fold composition of   derivations in the associative algebra of   differential operators.
This turns $\lieg$ into a \emph{restricted Lie algebra}.
The interplay between $p$-power map, Lie bracket and scalar multiplication is regulated by three axioms
(\cite[Chapter II, \S7]{Demazure-Gabriel1970}). In our situation, the $p$-power map on the Lie algebra
coincides with the Frobenius on cohomology, so the Bockstein kernel $\lieg_\red=H^1(Y,\O_Y)[\beta]$ is a restricted Lie subalgebra.
Here we are interested in
the \emph{Bockstein cokernel} $\lieh=\lieg/\lieg_\red$, with its  inherited  structure of   restricted Lie algebra.

Recall that for every algebraic group scheme $H$, the relative Frobenius is a homomorphism $F\colon H\ra H^{(p)}$
of group schemes, and its kernel is a  local group scheme $H[F]$. Local group schemes $H$
with $H=H[F]$ are called \emph{of height $\leq 1$}. 
According to \cite[Chapter II, \S7, Theorem 4.2]{Demazure-Gabriel1970}, the functor $H\mapsto \Lie(H)$ 
is  an equivalence between the categories of local groups schemes 
of height $\leq 1$ and finite-dimensional restricted Lie algebras.

\begin{proposition}
\mylabel{local pic multiplicative}
Suppose that the Frobenius map on $H^2(Y,\O_Y)$ has maximal Hasse--Witt rank.
Then the   group scheme $H$ of height $\leq 1$ corresponding to the Bockstein cokernel $\lieh=H^1(Y,\O_Y)/H^1(Y,\O_Y)[\beta]$
is multiplicative. If the reduced part of $G=\Pic^0_{Y/k}$ is geometrically reduced,
then the resulting local group scheme  $L=G/G_\red$ is also multiplicative. 
\end{proposition}

\proof
It suffices to treat the case that $k$ is perfect.
According to Mumford's analysis in \cite[Lecture 27]{Mumford1966},
the Lie algebra of the smooth connected group scheme $G_\red$ coincides with
the Bockstein kernel. In turn, the Bockstein cokernel is the Lie algebra for
the local group scheme $L=G/G_\red$, which gives an identification $H=L[F]$. But $L$ is multiplicative if and only
if its Frobenius kernel $H$ is multiplicative, because the higher Frobenius kernels $L[F^i]$ 
give a filtration on $L$ whose subquotients are isomorphic to $H$.

Seeking a contradiction, we assume that the inclusion $H^\mult\subset H$ is not an equality.
Since $k$ is perfect, the projection $H\ra H/H^\mult$ admits a section, and we get $H=H^\mult\oplus U$
for some non-trivial   unipotent local group scheme $U$.
According to \cite[Chapter IV, \S2, Corollary 2.13]{Demazure-Gabriel1970},
the Lie algebra $\Lie(U)$ contains a vector $b\neq 0$ with $b^{[p]}=0$.
Choose a vector  $a\neq 0$ from $ \lieg=H^1(Y,\O_Y)$ mapping to   $b$. Then $a^{[p]}\in \lieg_{\red}$.

Now recall that the $p$-power map equals the Frobenius map.
By construction $a$ does not lie in the Bockstein kernel.
Hence there is a largest integer $r\geq 0$  such that $a$ belongs to  the image of   $H^1(Y,W_r)\ra H^1(Y,\O_Y)$.
For the Bockstein operator, this means $\beta_r(a)\neq 0$.
Since $F(a)=a^{[p]}$ lies in the Bockstein kernel, we have $0=\beta_r(F(a))=F(\beta_r(a))$.

By definition, the range of the Bockstein $\beta_r$ is a quotient of $H^2(V^r_{r+1})$,
and $V^r_{r+1}$ consists of tuples $(0,\ldots,0,\lambda)$. In turn, there is
an identification $V^r_{r+1}=\O_Y$ of abelian sheaves, compatible with Frobenius.
By assumption, the Frobenius is bijective on $H^2(Y,\O_Y)$. 
With Lemma \ref{rank in sequence} below we infer that it is also bijective 
on the   the range  of the Bockstein operator. This gives $\beta_r(a)=0$, contradiction.
\qed

\medskip
In the above arguments, we have used the following simple observation:

\begin{lemma}
\mylabel{rank in sequence}
Let $V'\ra V\ra V''$ be an exact sequence of $k[F]$-modules whose underlying $k$-vector spaces are finite-dimensional.
If $V'$ and $V''$ have maximal Hasse--Witt rank, the same holds for $V$.
\end{lemma}

\proof
It suffices to treat the case that $k$ is perfect. We then have to check that the $p$-linear map $F:V\ra V$ is surjective.
Clearly, its image contains the image of $V'\ra V$.
This reduces us to the case $V'=0$, such that $V\subset V''$.
Choose a vector space basis $a_1,\ldots,a_r\in V$ and extend it to a basis $a_1,\ldots,a_n\in V''$.
Then $F(a_j)=\sum\lambda_{ij}a_i$ defines a matrix $A''=(\lambda_{ij})$. In turn, the Hasse--Witt determinant 
$\det(A'')$ is non-zero.
Since $V\subset V''$ is an $F$-invariant subspace, we see that $\lambda_{ij}=0$ for $1\leq j\leq r<i$,
and conclude that also $F:V\ra V$ has non-zero Hasse--Witt determinant.
\qed

\medskip
\noindent\emph{Proof of Theorem \ref{hasse--witt}.}
Suppose that the Frobenius has maximal Hasse--Witt rank   on $H^2(Y,\O_Y)$.
We have to check that $\Upsilon_G^0=0$ for the algebraic group scheme $G=\Pic^\tau_{Y/k}$.
In light of Proposition \ref{base-change vanish}, it suffices to treat the case that $k$ is perfect.
Then $G^0_\red$ is a subgroup scheme, and by Theorem  \ref{structure upsilon} we have to verify that
the local group scheme $L=G^0/G^0_\red$ is multiplicative.
This holds by Proposition \ref{local pic multiplicative}.
\qed

\section{The case of surfaces}
\mylabel{Surfaces}

We keep the assumptions of the previous section, so that $\Upsilon_{Y/k}$ is the maximal unipotent
quotient of $G=\Pic^\tau_{Y/k}$, where $Y$ is a proper scheme over our ground field $k$
of characteristic $p>0$. The goal now is to apply the general results of the previous section
to certain classes of $Y$, and establish   vanishing results.
For simplicity, we assume that $Y$ is equi-dimensional, of dimension $n\geq 0$.
Our first observation is:

\begin{proposition}
\mylabel{vanishing curves}
We have $\Upsilon_{Y/k}=0$ provided that $Y$ is a curve.
\end{proposition}
 
\proof 
The N\'eron--Severi group $\NS(X\otimes k^\alg)$ is a free group, according to  
\cite[Section 9.4, Corollary 14]{Bosch-Luetkebohmert-Raynaud1990}.
Hence by Proposition \ref{galois module upsilon}, the group scheme of components for $\Upsilon_{Y/k}$ is trivial.
Furthermore, the group $H^2(Y,\O_Y)$ vanishes by dimension reason, so the Frobenius has a priori
maximal Hasse--Witt rank. The  the assertion thus follows 
from Theorem \ref{hasse--witt}.
\qed

\medskip
\newcommand{\trace}{\operatorname{tr}}
Let $\omega_Y$ be the \emph{dualizing sheaf}, with its trace map $\trace\colon H^n(Y,\omega_Y)\ra k$.
In turn,  for every coherent sheaf $\shF$ we get a pairing
$$
H^i(Y,\shF)\times\Ext^{n-i}(\shF,\omega_Y)\stackrel{\can}{\lra} H^n(Y,\omega_Y)\stackrel{\trace}{\lra} k,
$$
which is non-degenerate for $i=n$,   regardless of the singularities.
The pairings remain non-degenerate in all degrees $i\leq n$ provided that $Y$ 
is Cohen--Macaulay. We record:

\begin{proposition}
\mylabel{vanishing duality}
We have $\Upsilon_{Y/k}^0=0$ provided that  the scheme  $Y$ has dimension $n\geq 2$, is Cohen--Macaulay,
and has the property  $H^{n-2}(Y,\omega_Y)=0$.
\end{proposition}
 
\proof
The above Serre Duality gives $h^2(\O_Y)=h^{n-2}(\omega_Y)=0$, and the assertion follows
from Theorem \ref{hasse--witt}.
\qed

\medskip
Now suppose that $Y$ is a surface. For every invertible sheaf $\shL$,
the Euler characteristic $\chi(\shL^{\otimes t})$ is a numerical polynomial   of degree two in the variable $t$,
which can be written as 
$$
\chi(\shL^{\otimes t}) = \frac{(\shL\cdot\shL)}{2} t^2 - \frac{(\shL\cdot \omega_Y)}{2} t + \chi(\O_Y).
$$
The quadratic term is determined by  the self-intersection number $(\shL\cdot\shL)$, whereas  
\emph{the linear term   defines an  integer} $(\shL\cdot\omega_Y)$,
which coincides with the usual intersection number of invertible sheaves provided that  $Y$ is Gorenstein.  

\begin{theorem}
\mylabel{vanishing surfaces}
Let  $Y$ be a reduced surface.
Suppose there is an invertible sheaf $\shL$ such that for each irreducible component $Z\subset Y$,
the restriction $\shL_Z$ is nef and has   $(\shL_Z\cdot\omega_Z)<0$.
Then the cohomology group $H^2(Y,\O_Y)$ vanishes, the N\'eron--Severi group $\NS(Y\otimes k^\alg)$ is free,
and  the maximal unipotent finite quotient  $\Upsilon_{Y/k}$ is trivial.
\end{theorem}

\proof 
\newcommand{\uHom}{\underline{\Hom}}
The statement on the maximal unipotent quotient follows
from the assertions  on   cohomology and   N\'eron--Severi group, using Proposition \ref{galois module upsilon} and Theorem \ref{hasse--witt}.
To proceed, we first reduce to the case   of integral   surfaces that are Cohen--Macaulay.
Let $X_1,\ldots,X_r$ be the $S_2$-ization of the irreducible components $Y_1,\ldots,Y_r\subset Y$,
and write $X$ for their disjoint union (see for example \cite{Schroeer-Vezzosi2004} for details on the $S_2$-ization). 
The resulting  finite morphism $f\colon X\ra Y$ is surjective, hence
induces an injection on N\'eron--Severi groups. Moreover, the cokernel in
the   short exact sequence $0\ra \O_Y\ra f_*(\O_X)\ra\shF\ra 0$ is at most zero-dimensional 
(see \cite[corollaire~5.10.15]{EGAIVb}).
In the resulting long exact sequence
$$
H^1(Y,\shF)\lra H^2(Y,\O_Y)\lra H^2(X,\O_X)\lra H^2(Y,\shF),
$$
the outer terms vanish, and we get an identification  $H^2(Y,\O_Y)=H^2(X,\O_X)$. 
Now fix some irreducible component $Z=Y_i$, and consider the resulting connected
component $Z'=X_i$.  
Its dualizing sheaf  is given by $f_*(\omega_{Z'})=\uHom(f_*(\O_{Z'}),\omega_Z)$,
and we get a short exact sequence $0\ra f_*(\omega_{Z'})\ra\omega_Z\ra\shG\ra 0$ where the
cokernel is at most zero-dimensional. Tensoring with $\shL^{\otimes t}$ and comparing   linear terms in
the numerical polynomials, we see that $(\shL_{Z'}\cdot \omega_{Z'})=(\shL_Z\cdot\omega_Z)$.
Clearly, the pullback of $\shL$ under the finite surjection $Z'\ra Z$ remains  nef.
Summing up, it suffices to treat the case that $Y$ is integral and Cohen--Macaulay. 

Next, we verify that $h^2(\O_Y)=0$. Seeking a contradiction,  we assume that there is a short exact sequence
$0\ra\O_Y\ra\omega_Y\ra\shF\ra0$, where $\shF$ is a torsion sheaf.
Tensoring with $\shL^{\otimes t}$ and using Serre Duality, we get
\begin{equation}
\label{euler characteristics}
\chi(\shL^{\otimes t}) + \chi(\shL^{\otimes t}\otimes\shF) = \chi(\shL^{\otimes t}\otimes\omega_Y)= \chi(\shL^{\otimes -t}).
\end{equation}
Let $\zeta_i\in Y$ be the generic points of $\Supp(\shF)$
whose closures $C_i\subset Y$ are one-dimensional, write $m_i\geq 1$ for the length of the finite $\O_{Y,\zeta_i}$-module
$\shF_{a_i}$, and let $C$ be the union $\bigcup C_i$. By the results in  \cite[Section II.2]{Kleiman1966},  
we have $\chi(\shL^{\otimes t}\otimes\shF) = \sum m_i (\shL\cdot C_i)t+\chi(\O_C)$, with $\sum m_i (\shL\cdot C_i)\geq 0$.
In equation \eqref{euler characteristics}, the linear term  on the left has coefficient 
$-\frac{1}{2}(\shL\cdot\omega_Y) +  \sum m_i (\shL\cdot C_i)>0$,
whereas the  coefficient on the right is $\frac{1}{2}(\shL\cdot \omega_Y)<0$, contradiction.
This shows $h^2(\O_Y)=0$.

It remains to verify the assertion on the N\'eron--Severi group. Let $f\colon X\ra Y$ be some resolution of singularities.
This proper surjective morphism induces an inclusion  $\NS(Y)\subset\NS(X)$. Furthermore, we obtain a short exact sequence
$0\ra f_*(\omega_X) \ra \omega_Y\ra\shF\ra 0$ for some torsion sheaf $\shF$. Tensoring with $\shL^{\otimes t}$,
and using Serre duality alongside the Leray--Serre spectral sequence,  we get
$$
\chi(\shL^{\otimes -t}) = \chi(\shL_X^{\otimes -t})-\chi(\shL^{\otimes t}\otimes R^1f_*\O_X) + \chi(\shL^{\otimes t}\otimes\shF).
$$
Looking at the linear terms and using the notation of the preceding paragraph, 
we get
$$\frac{1}{2}(\shL\cdot \omega_Y) = \frac{1}{2} (\shL_X\cdot\omega_X) + \sum m_i (\shL\cdot C_i),$$
and conclude that
$(\shL_X\cdot\omega_X)<0$. This reduces us to the case of regular irreducible surfaces $Y$.
Likewise, one easily reduces to the case   that $k$ is separably closed.

The base-change to the algebraic closure $k^\alg$ is not necessarily regular or normal, not even reduced.
Let $X\ra Y\otimes k^\alg$ be the normalization of the reduction, and consider the composite morphism
$f\colon X\ra Y$. According to a result of Tanaka (\cite[Theorem~4.2]{Tanaka2018}, see also \cite[Theorem 1.1]{Patakfalvi-Waldron2017}),
we have the following equality $\omega_X=f^*(\omega_Y)\otimes \O_X(-R)$ for some curve $R\subset X$.
In turn, $(\shL_X\cdot\omega_X)<0$. Now let $r\colon S\ra X$ be the minimal resolution of singularities,
with exceptional divisor $E=E_1+\ldots+E_r$, and $K_{S/X}=\sum\lambda_iE_i$ be the unique $\QQ$-divisor
with $ (K_{S/X}\cdot E_i)= (K_S\cdot E_i)$. Since the resolution is minimal, we must have $\lambda_i\leq 0$.
Now choose some Weil divisor  $K_X$ representing $\omega_X$, and consider the rational pullback $r^*(K_X)$ in the
sense of Mumford \cite{Mumford1961}, compare also \cite{Schroeer2019}. We then have
$$
(\shL_S\cdot \omega_S) = f^*(\shL)  \cdot  ( K_{S/X}+ f^*K_X )= (\shL_X\cdot \omega_X)<0.
$$
Using that $\shL_S$ is nef, we infer that the plurigenera $ h^0(\omega_S^{\otimes n})$, $n\geq 1$
of the smooth surface $S$ vanish, so its Kodaira dimension must be $\operatorname{kod}(S)=-\infty$.
By the Enriques classification, the surface is either $S=\PP^2$ or admits a ruling.
In both cases,  $\NS(S)$ is free. In turn, the same holds for the subgroup $\NS(Y\otimes k^\alg)$.
\qed

\medskip
This applies in particular for reduced \emph{del Pezzo surfaces} $Y$, which by definition 
are Gorenstein, with $\omega_Y$ anti-ample, and have   $h^0(\O_Y)=1$.

\begin{corollary}
\mylabel{del Pezzo corollary}
Let $Y$ be a reduced del Pezzo surface. Then $\Upsilon_{Y/k}$ is trivial.
\end{corollary}

Note that Miles Reid \cite{Reid1994} has classified reduced non-normal del Pezzo surfaces over  algebraically closed fields.
Then the algebraic group
$G=\Pic^\tau_{Y/k}$ is smooth, and there are cases with $H^1(Y,\O_Y)\neq 0$.
 
Tanaka constructed Mori fiber spaces in characteristic $p\leq 3$, where the generic fiber 
is a normal projective  surfaces $Y$ with $h^0(\O_Y)=1$ having
only $\QQ$-factorial klt-terminal singularities, the    $\QQ$-divisor $K_Y$ is anti-ample,
yet the Picard group contains elements  of order $p$  (\cite[Theorem 1.2]{Tanaka2016}, with further investigation 
in \cite{Bernasconi-Tanaka2019}).
From Theorem \ref{vanishing surfaces}, we see that such elements  must come from the   component $\Pic^0_{Y/k}$ of the origin.

Note also that the reducedness assumption in our results is indispensable:
Suppose that $S$ is an irreducible  smooth  surface, and let $\shL$ be an invertible sheaf
such that its dual is ample. After passing to suitable multiples, we achieve 
that $h^1(\shL^\vee)=h^2(\shL^\vee)=0$, and 
$\shL \otimes\omega _S$ becomes anti-ample.
Consider the quasicoherent $\O_S$-algebra $\shA=\O_Y\oplus\shL^\vee$,
with multiplication $(f,s)\cdot (f',s') = (ff', fs'+f's)$, and let $Y=\Spec(\shA)$
be its relative spectrum. The structure morphism $f\colon Y\ra S$ has a canonical section,
which is given by the projection $\shA\ra \O_Y$ and identifies $S$ with $Y_\red$.
One also says that $Y$ is a \emph{ribbon} on $S$.
The resulting short exact sequence of abelian sheaves $0\ra \shL^\vee\ra\O_Y^\times\ra\O_S^\times\ra 1$
shows that the inclusion $S\subset Y$ induces an identification of Picard schemes.
Moreover, the relative dualizing sheaf for $f\colon Y\ra S$ is given by
$f_*(\omega_{Y/S})=\uHom(\shA,\O_Y)=\shA\otimes\shL$.
In turn, we have $\omega_Y=f^*(\shL\otimes\omega_S)$, which is anti-ample.
Summing up,
$Y$ is an irreducible  non-reduced del Pezzo surface with $\Upsilon_Y=\Upsilon_S$.
The latter easily becomes non-trivial, e.g.~if $S$ 
is a simply-connected Enriques surface in characteristic $p=2$.

\section{Enriques surfaces and del Pezzo surfaces}
\mylabel{Enriques del Pezzo}

\newcommand{\hS}{{S'}}
\newcommand{\tS}{{\tilde{S}}}
\newcommand{\tZ}{{\tilde{Z}}}
\newcommand{\hD}{{\hat{D}}}

Let $k$ be an algebraically closed ground field of characteristic $p=2$.
In this section, we construct certain  normal  Enriques  surfaces  $Z$    where the numerically trivial part of the Picard scheme
is unipotent of  order two,
together with a finite universal homeomorphism $\nu\colon Z'\ra Z$ of degree two from  a normal del Pezzo surface 
of degree four with Picard number one. Theses surfaces will arise from very special simply-connected Enriques surfaces.

Throughout, $S$ denotes an \emph{Enriques surface}.
This means that $S$ is a regular connected surface with $c_1=0$ and $b_2=10$.
The group $\Num(S)=\Pic(S)/\Pic^\tau(S)$ of numerical classes is  a free abelian group of rank $\rho=10$
called the \emph{Enriques lattice}. The intersection form is isomorphic to $E_8\oplus H$,
comprising the root lattice of type $E_8$  and the hyperbolic lattice $H=(\begin{smallmatrix}0&1\\1&0\end{smallmatrix})$.
The group scheme $P=\Pic^\tau_{S/k}$ of numerically trivial invertible sheaves is finite
of order two. In characteristic $p=2$, there are three possibilities for $P$, namely 
$$
\mu_2\quadand \ZZ/2\ZZ\quadand \alpha_2.
$$
The respective Enriques surfaces $S$ are called \emph{ordinary}, \emph{classical} and \emph{supersingular}.
The inclusion $P\subset\Pic_{S/k}$ yields a $G$-torsor $\epsilon\colon \tS\ra S$, 
where $G=\underline{\Hom}(P,\GG_m)$
is the Cartier dual. For ordinary $S$   the Cartier dual  $G=\ZZ/2\ZZ$
is \'etale, and the total space $\tS$ is a \emph{K3 surface}.
In the other two cases  $P$ is unipotent,
 $G$ is local,  and the integral Gorenstein surface $\tilde{S}$ necessarily acquires singularities.
One   says that $S$ is a \emph{simply-connected Enriques surface}, and that $\tS$ is its \emph{K3-like covering}.
Note that the relation between inclusions $P\subset \Pic_{S/k}$
and $G$-torsors $\tS\ra S$ goes back to  Raynaud \cite{Raynaud1970}.
The trichotomy of Enriques surfaces for $p=2$ was developed  by Bombieri and Mumford \cite{Bombieri-Mumford1976}.
For more information on K3-like coverings, see for example \cite{Schroeer2017}.

We write $\gamma\colon S'\ra\tilde{S}$ for the normalization. Since $S$ is regular and $S'$ and $\tS$ are Cohen--Macaulay,
both projections   $\nu\colon S'\ra S$ and $\epsilon\colon \tS\ra S$ are finite and flat of degree two (\cite[proposition~6.1.5]{EGAIVb}).
Hence all fibers take the form  $\Spec\kappa(s)[\epsilon]$, which ensures that both  $S'$ and $\tS$ are Gorenstein.
The relative dualizing sheaf for $\gamma\colon S'\ra\tilde{S}$ is given by 
$$
\gamma_*(\omega_{S'/\tS})=\underline{\Hom}(\gamma_*(\O_{S'}),\O_\tS).
$$
It is also the  ideal sheaf for the \emph{ramification locus} $R\subset S'$ of the normalization map. This Weil divisor must be Cartier, because $S'$ and $\tS$
are Gorenstein.
As explained in \cite{Ekedahl-Shepherd-Barron2004}, there is a unique effective divisor  $C\subset S$ with $R=\nu^{-1}(C)$.
Ekedahl and Shepherd-Barron call $C\subset S$ the \emph{conductrix} of the  Enriques surface
$S$. Note that   $S$ is    simply-connected   if $C$ is non-empty.

Now suppose that 
$E=E_1+\ldots+E_r$ is a configuration of $(-2)$-curves whose intersection matrix $(E_i\cdot E_j)$ is negative-definite.
We also say that $E$ is an \emph{ADE-curve}. 
Let $f\colon S\ra Z$ be its contraction. The resulting surfaces $Z$ are called \emph{normal Enriques surfaces}.
They are projective and  their  singularities are rational double points, such that   
$\omega_Z=f_*(\omega_S)$ and $\omega_S=f^*(\omega_Z)$.
By the Hodge Index Theorem, we have $r\leq 9$.

\begin{proposition}
\mylabel{numerical group}
The numerical group $\Num(Z)$ if free of rank $\rho=10-r$,
and the homomorphism of group schemes $f^*\colon \Pic^\tau_{Z/k}\ra\Pic^\tau_{S/k}$ is an isomorphism.
\end{proposition}

\proof
The Picard group of $Z$ can be viewed as the orthogonal complement
of the curves $E_1,\ldots,E_r\in \Pic(S)$, via the preimage map $f^*\colon\Pic(Z)\ra \Pic(S)$.
It follows that $\Pic(Z)$ is finitely generated of  rank $\rho=10-r$, and the assertion on the numerical group $\Num(Z)$ follows.

The subgroup $\Pic^\tau(Z)\subset\Pic(Z)$ is the torsion part. This is cyclic of order two 
provided that $S$ is classical, and then the assertion on $\Pic^\tau_{Z/k}$ is already contained in  the preceding paragraph.
It remains to treat the case that $S$ is ordinary or supersingular. In other words, the group scheme $P=\Pic^\tau_{Z/k}$ is local
of height one, with one-dimensional tangent space $\Lie(P)=H^1(S,\O_S)$. 
The Leray--Serre spectral sequence for the contraction $f\colon S\ra Z$ gives an exact sequence
$$
0\lra H^1(Z,\O_Z)\lra H^1(S,\O_S)\lra H^0(Z,R^1f_*(\O_S)).
$$
The term on the right vanishes, because $Z$ has only rational singularities, and it follows
that $\Pic_{Z/k}\ra \Pic_{S/k}$ induces a bijection on tangent spaces.
This homomorphism of groups schemes must be a monomorphism, because $\O_S=f_*(\O_X)$.
In turn, the inclusion $\Pic^\tau_{Z/k}\subset\Pic^\tau_{S/k}$ is an isomorphism.
\qed

\medskip
Set $P=\Pic^\tau_{S/k}=\Pic^\tau_{Z/k}$ and let $G=\underline{\Hom}(P,\GG_m)$ be the Cartier dual.
The inclusion of $P\subset\Pic_{Z/k}$  
corresponds to   a $G$-torsor $\tZ\ra Z$, and the total space $ Z'$ is an integral Gorenstein surface.
Let $ Z'\ra \tZ$ be its normalization. Of course, we have analogous constructions for the
Enriques surface $S$, and by naturality we get a commutative diagram
\begin{equation}
\label{normalization and k3-like}
\begin{CD}
S'		@>>>	\tS	@>>> 	S\\
@VVV		@VVV		@VVf V\\
 Z'		@>>>	\tZ	@>>>	Z
\end{CD}
\end{equation}
of projective integral surfaces. We write $\nu\colon Z'\ra Z$ for the composition of the lower arrows.
The image $D=f(C)$   is a Weil divisor on 
the normal surface $Z$. We call it the \emph{conductrix of the normal Enriques surface}.

\begin{theorem}
\mylabel{del pezzo}
Suppose   the Enriques surface $S$ and the chosen ADE-curve $E\subset S$ satisfy  the following three conditions:
\begin{enumerate}
\item
The conductrix $C\subset S$ is not supported by the exceptional divisor $E$.
\item
There are integers $m_1,\ldots,m_r$ such that $(C\cdot E_j)=\sum m_i(E_i\cdot E_j)$ for all $1\leq j\leq r$.
\item
The ADE curve $E=E_1+\ldots+E_r$ has $r=9$ irreducible components.
\end{enumerate}
Then  the conductrix $D\subset Z$ is Cartier, and   $Z'$ is a normal del Pezzo surface of degree 
$K_{Z'}^2=2D^2$, with  canonical class $K_{Z'}=-\nu^{-1}(D)$, irregularity $h^1(\O_{Z'})=0$ and Picard group $\Pic(Z')=\ZZ$.
Moreover, all singularities on $ Z'$ are rational double points.
\end{theorem}

\proof
Since the horizontal maps in the diagram \eqref{normalization and k3-like} are universal homeomorphisms
and the scheme $Z$ is $\QQ$-factorial,
the schemes $ Z',\tZ$ and $Z$ have the same Picard number, so Condition (iii) yields $\rho( Z')=\rho(Z)=1$.
Condition (i) ensures that the Weil divisor $D\subset Z$ is non-empty, and (ii) means that it is Cartier.

By definition, the preimage of the  conductrix $C\subset S$ on $S'$ is the ramification locus for the normalization $S'\ra \tS$.
Let $U\subset\tZ$ be the preimage of the regular locus $\Reg(Z)$.
Since $\tS\ra\tZ$ becomes an isomorphism over  $U$, 
the preimage of $D$ in $\tZ$ and the branch locus in $\tZ$  for the normalization $ Z'\ra\tZ$ coincide,  at least over $U$.
But the preimage $\nu^{-1}(D)$ and the ramification curve $R\subset Z'$ have no embedded components, and $Z'$ is normal, so 
$\nu^{-1}(D)=R$. As $D\subset Z$ is Cartier, the same holds for $R\subset Z'$.
Since $\omega_\tZ$ and $\omega_{ Z'/\tZ}=\O_{ Z'}(-R)$ are invertible, the normal surface $ Z'$ is Gorenstein.
We have  $-(K_{Z'}\cdot R) = \nu^*(D)^2 = \deg(\nu)\cdot D^2=2D^2>0$.
By the Nakai criterion, $\omega_{Z'}$ is anti-ample, thus $ Z'$ is a normal del Pezzo surface.

Seeking a contradiction, we now assume that $H^1( Z',\O_{Z'})\neq 0$.
This is the tangent space to the Picard scheme, so $A=\Pic^0_{ Z'/k}$ is non-zero.
The latter   is smooth, because the obstruction group $H^2( Z',\O_{Z'})\simeq H^0( Z',\omega_{Z'})$
vanishes, compare \cite[Lecture 27]{Mumford1966}.
In turn, the group scheme $A\neq 0$ is an abelian variety, and we conclude that for each prime $l\neq 2$,
there is a $\mu_l$-torsor $ Z''\ra Z'$ with   connected total space.
The projection $\nu\colon Z'\ra Z$ is a universal homeomorphism, so by \cite[expos\'e~IX, th\'eor\`eme~4.10]{SGA1}, 
the finite \'etale Galois covering $ Z''\ra  Z'$ is the base-change of some finite \'etale Galois covering of $Z$.
This implies that $\Pic(Z)$ contains an element of order $l$, contradiction. Thus $h^1(\O_{Z'})=0$.

It then follows $\Pic( Z')=\ZZ$ by  \cite[Lemma 2.1]{Schroeer2001}.  
It remains to check that the singularities on $ Z'$ are  rational double points.
Since $ Z'$ is Gorenstein, the task is to show that they are rational.
Seeking a contradiction, we assume that there is at least one non-rational singularity.
Consider the minimal resolution of non-rational singularities
$r\colon Y\ra Z'$. According to \emph{loc.~cit.}, Theorem 2.2, there exists  a fibration $\varphi\colon Y\ra B$ over some curve
of genus $g>0$, and $r\colon Y\ra  Z'$ is the contraction of some section $E\subset Y$.
Set $V=\Reg(Z)$, and fix some prime $l\neq 2$ such that the finitely generated abelian group $\Pic(V)=\Pic(S)/\sum\ZZ E_i$
contains no element of order $l$.
As in the preceding paragraph, we find some $\mu_l$-torsor $Y'\ra Y$, 
which yields a finite \'etale Galois covering of $V$. This implies that $\Pic(V)$
contains an element of order $l$, contradiction.
\qed

\begin{corollary}
\mylabel{non-regular}
Keep the assumptions of Theorem~\ref{del pezzo}. Then the normal del Pezzo surface $ Z'$ is not regular, we have $\Pic(Z')=\ZZ K_{Z'}$,
and the Cartier divisor $D\subset Z$ has selfintersection  $1\leq D^2\leq 4$.
\end{corollary}

\proof
Suppose $ Z'$ were regular. It contains no $(-1)$-curves,
because $\rho=1$. According to \cite[Theorem 8.1.5]{Dolgachev2012}, the regular del Pezzo surface
$ Z'$ is isomorphic to $\PP^2$. But $K^2_{\PP^2}=9$ is odd,
whereas $K_{ Z'}^2=2D^2$ is even, contradiction.

Let $r\colon X\ra  Z'$ be the minimal resolution of singularities.
Then $X$ is a \emph{weak del Pezzo surface},
which for integral surfaces means that it is Gorenstein and the inverse   of the dualizing sheaf 
is  nef and big.
Let $Y$ be a minimal model, obtained from a successive contraction $X=X_0\ra\ldots\ra X_n=Y$   of $(-1)$-curves.
Then $Y$ is either the projective plane or a Hirzebruch surface, and in both cases we have $K_Y^2\leq 9$.
This gives  $2D^2=K_{Z'}^2=K_X^2=K_Y^2-n\leq 9$, and in turn $1\leq D^2\leq 4$.
In particular,  the degree $K_{Z'}^2$ of the normal del Pezzo surface belongs to the set $\{2,4,6,8\}$.

Suppose the canonical class does not generate the Picard group. Then we are in the situation
$K_{Z'}=2A$ for some Cartier divisor $A$, and either $A^2=1$ or $A^2=2$.
In case $K_{Z'}^2=8$, the surface $X$ is the Hirzebruch surface with invariant $e=2$,
and $Z'$ is obtained by contraction of the $(-2)$-curve. It then follows that $K_{Z'}$ generates
the Picard group, contradiction.
Now suppose that $K_{Z'}^2=4$, with 
intersection numbers $A^2=1$ and $A\cdot K_{Z'}=-2$.
In turn, the Euler characteristic  $\chi(\shL)=A\cdot(A-K_{Z'})/2+\chi(\O_{Z'})$  
of the invertible sheaf $\shL=\O_{Z'}(A)$ is not an integer, contradiction.
\qed

\medskip
For later use, we also record the following vanishing result:

\begin{corollary}
\mylabel{raynaud vanishing}
Keep the assumptions of Theorem~\ref{del pezzo}. Then $h^1(\omega_{Z'}^{\otimes -i})=0$ 
for all  integers $i$.
\end{corollary}

\proof
By Serre Duality, it suffices to verify this for $i\geq 0$.
Let $r\colon X\ra  Z'$ be the minimal resolution of singularities. Since $ Z'$ has only rational double points,
we have $\omega_X^{\otimes -i}=r^*(\omega_{Z'}^{\otimes- i})$ and $h^1(\omega_X^{\otimes-i})=h^1(\omega_{Z'}^{\otimes -i})$.
The regular surface $X$ is obtained from its minimal model $Y$, which is the projective plane or a Hirzebruch surface,
by a sequence of blowing-ups of closed points. It follows that $X$ lifts to characteristic zero.
In particular, we may apply Raynaud's vanishing result (\cite[Corollary 2.8]{Deligne-Illusie1987})
$$
H^1(X,\omega_X^{\otimes -i} ) = H^1(X,\Omega_{X/k}^2\otimes \shL) =0
$$
 for the nef and big invertible sheaf $\shL=\omega_X^{\otimes-(i+1)}$.
\qed

\section{Exceptional Enriques surfaces}
\mylabel{Exceptional enriques}

Now suppose that $S$ is a simply-connected Enriques surface   
whose conductrix takes the form
$$
C=2C_1+3C_3+5C_4+2C_2+4C_5+4C_6+3C_7+3C_8+2C_0+C_9,
$$
where the ten irreducible components $C_0,\ldots, C_9$ are  $(-2)$-curves with simple normal crossings having the following   dual graph $\Gamma$:
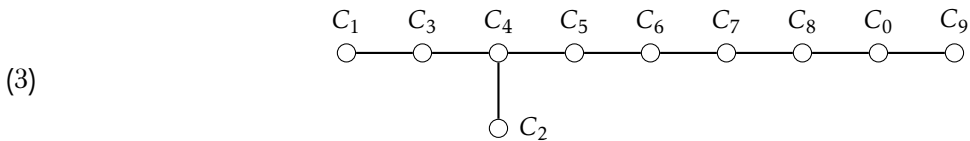
\begin{equation}
\label{dual graph 237}
\begin{gathered}
\begin{tikzpicture}
[node distance=1cm, font=\small]
\tikzstyle{vertex}=[circle, draw, inner sep=0mm, minimum size=1.5ex]
\node[vertex]	(C1)  	at (0,0) 	[label=above:{$C_1$}] 			{};
\node[vertex]	(C3)			[right of=C1, label=above:{$C_3$}]	{};
\node[vertex]	(C4)			[right of=C3, label=above:{$C_4$}]	{};
\node[vertex]	(C2)			[below of=C4, label=right:{$C_2$}]	{};
\node[vertex]	(C5)			[right of=C4, label=above:{$C_5$}]	{};
\node[vertex]	(C6)			[right of=C5, label=above:{$C_6$}]	{};
\node[vertex]	(C7)			[right of=C6, label=above:{$C_7$}]	{};
\node[vertex]	(C8)			[right of=C7, label=above:{$C_8$}]	{};
\node[vertex]	(C0)			[right of=C8, label=above:{$C_0$}]	{};
\node[vertex]	(C9)			[right of=C0, label=above:{$C_9$}]	{};

\draw [thick] (C1)--(C3)--(C4)--(C5)--(C6)--(C7)--(C8)--(C0)--(C9);
\draw [thick] (C4)--(C2);
\end{tikzpicture}
\end{gathered}
\end{equation}
One also says that $S$ is a \emph{exceptional Enriques surface  of type $T_{2,3,7}$}.
Here the indices $2,3,7$ denote the length of the terminal chains in the star-shaped tree $\Gamma$, including the
central vertex $C_4\in\Gamma$.
Such  Enriques surfaces were already considered  in the monograph of Cossec and Dolgachev \cite[Chapter III, \S4]{Cossec-Dolgachev1989}.

The general notion of \emph{exceptional Enriques surfaces} was introduced by 
Ekedahl and Shepherd-Barron \cite{Ekedahl-Shepherd-Barron2004}, who  
studied them in detail. They can be characterized in terms of  the conductrix $C\subset S$,
and also by properties of the Hodge ring $\bigoplus_{ij} H^i(S,\Omega_S^j)$.
Explicit equations for birational models   were found by Salomonsson \cite{Salomonsson2003}.
For examples,  the equation
$$
z^2+(y^4+x^4)x^3y^3s^4 + \lambda x^5y^3s^3t+xyt =0,\quad\lambda\neq 0,
$$
as well as 
$$
z^2+x^3y^7s^4+\mu x^8s^3t + xyt^4=0,\quad \mu\neq 0
$$
define  birational models for  
exceptional Enriques surfaces of type $T_{2,3,7}$ as inseparable double covering of a Hirzebruch surface  with coordinates $x,y,s,t$.
The first equation gives classical, the second equation supersingular  Enriques surfaces.

The reduced curve $F=C_0+\ldots+C_8$ on the Enriques surface $S$ supports a curve of canonical type with Kodaira symbol $\II^*$. 
Let $\varphi\colon S\ra\PP^1$ be the resulting
genus-one fibrations. This fibration is quasielliptic, and there is no other genus-one fibration,
according to  \cite[Theorem C]{Ekedahl-Shepherd-Barron2004}. The fiber corresponding to $F$ 
is  multiple, because otherwise  $2\leq (C_9\cdot F)= (C_9\cdot C_0)=1$, contradiction.
Since $b_2=10$, all other fibers are irreducible,  thus have Kodaira symbol $\II$. If $S$ is classical, there must be another multiple fiber.
In the supersingular case, all other fibers are simple
(\cite[Theorem 5.7.2]{Cossec-Dolgachev1989}).
 
Let $f\colon S\ra Z$ be the contraction of the ADE-curves $C_1+\ldots+C_8$ and $C_9$.
Then $Z$ is a  normal  Enriques surface with $\Sing(Z)=\{a,b\}$, 
where the first local ring  $\O_{Z,a}$ is a rational double point of type
$E_8$, and the second local ring  $\O_{Z,b}$ is a rational double point of type $A_1$.
Write $D_0=f(C_0)$ for the image of the remaining $(-2)$-curve, which is a Weil divisor.
The conductrix of the normal Enriques surface is 
$D=f(C)=2D_0$.
One easily sees that 
Theorem \ref{del pezzo} applies, so $D\subset Z$ is Cartier, and we get a a normal del Pezzo surface $Z'$
  as  an inseparable double covering $\nu\colon Z'\ra Z$. The goal of this section is
to study the geometry of these surfaces in detail.

Following  Hartshorne \cite{Hartshorne1994}, we write 
$\APic(X)$ for the group
of isomorphism classes of reflexive rank-one sheaves, on a given normal noetherian scheme $X$.
It is called the  \emph{almost Picard group}, and could also be seen as the group of
1-cycles modulo linear equivalence. 
If $X$ is a proper surface, the group $\APic(X)$ is endowed with Mumford's \emph{rational selfintersection numbers}
\cite{Mumford1961}, which extend the usual intersection numbers for invertible sheaves. 
If $X=\Spec(R)$ is local, we   use the more traditional $\Cl(R)=\APic(X)$.
In the global case, we prefer  $\APic(X)$, because it emphasizes the relation to the Picard group.

\begin{proposition}
\mylabel{picard group downstairs}
The group $\APic(Z)$ is generated by $D_0$ and the canonical class $K_Z$,  with selfintersection number  $D_0^2=1/2$.
Moreover, the subgroup $\Pic(Z)$ has index two, and is generated by the conductrix   $D$, which has $D^2=2$,
together with  $K_Z$.
\end{proposition}

\proof
Via the pullback map $f^*\colon \Pic(Z)\ra\Pic(S)$, we 
may regard $\Pic(Z)$ as    the orthogonal complement of the nine curves $C_1,\ldots,C_9\in \Pic(S)$.
This orthogonal complement has rank one, and  is generated by the numerically trivial $K_Z$  and the linear combination
$$
2C_0 + 2(2C_8+3C_7+4C_6+5C_5+6C_4+3C_2+2C_1+4C_3) + C_9.
$$
This coincides with $f^*(D)$, and gives the selfintersection $D^2=f^*(D)^2=2$.
The almost Picard group can be seen as the cokernel for the inclusion $\sum_{i=1}^9\ZZ C_i\subset\Pic(S)$.
Since $C_0,\ldots,C_9\in\Pic(S)$ form a basis modulo the torsion part, the assertion on $\APic(Z)$ follows.
\qed

\medskip
Let $\tZ\ra Z$ be the canonical covering and $Z'\ra\tZ$ be its normalization, 
as considered  in the previous section.  Consider the composite morphism $\nu\colon Z'\ra Z$, which is finite of degree two.
According to Theorem \ref{del pezzo}, the total space $Z'$ is a normal del Pezzo surface of degree $K_{Z'}^2=4$
with $\Pic(Z')=\ZZ$.
In the next sections, we will embed $Z'$ into some normal threefold, and use $\nu$ as a gluing map for a denormalization.
Our goal here is to understand the geometry of the double covering $\nu\colon Z'\ra Z$.
The main task is to understand what happens over the conductrix   $D=f(C)$.

The following terminology will be useful: The \emph{rational cuspidal curve} is  the projective scheme
$\Spec k[t^2,t^3]\cup\Spec k[t^{-1}]$,
which is the integral singular curve of genus one whose local rings are unibranch.
A \emph{ribbon} on a scheme $X$ is a closed embedding $X\subset Y$
whose ideal sheaf $\shN\subset\O_Y$ satisfies $\shN^2=0$, such that the $\O_Y$-module $\shN$
is actually an $\O_X$-module, and that $\shN$ is invertible as $\O_X$-module.
This     terminology is due to Bayer and Eisenbud \cite{Bayer-Eisenbud1995}.

\begin{proposition}
\mylabel{ribbon downstairs}
The scheme $ D_\red$ is the rational cuspidal curve, and the conductrix $D=f(C)$ is a ribbon on $D_\red$ with ideal sheaf $\shN$
of degree $\deg(\shN)=-1$.
The cuspidal point of   $D_\red$ is located at the  $E_8$-singularity $a\in Z$.
\end{proposition}

\proof
\newcommand{\bS}{\bar{S}}
\newcommand{\ba}{\bar{a}}
\newcommand{\bD}{\bar{D}}
The morphism $f\colon S\ra Z$ factors over the contraction $g\colon S\ra \bS$ of the ADE-curve $C_1+\ldots+C_8$.
This creates a rational double point $\ba\in \bS$ of type $E_8$, and we have $\O_{\bS,\ba}=\O_{Z,a}$.
Since the normal surface $\bS$ is locally factorial, the integral curve   $\bD_0=g(C_0)$ remains Cartier.
It contains $\ba$, and the local ring $\O_{\bS,\ba}$ is singular, whence $\O_{\bD_0,\ba}$ is singular as well.
This singularity on $\bD_0$ must be unibranch, in light of the dual graph \eqref{dual graph 237}.
 
Since the   fibration $\varphi\colon S\ra\PP^1$ factors over $\bS$, we must have $(\bD_0\cdot \bD_0)=0$.
The Adjunction Formula gives $\deg(\omega_{\bD_0}) = (K_{\bS}+\bD_0)\cdot \bD_0 = 0$.
We have  $h^0(\O_{\bD_0})=1$, because the curve is integral, hence   $h^1(\O_{\bD_0})=1$.
The classification of integral curves of genus one shows that $\bD_0$ is the rational cuspidal curve.

The induced morphism $h\colon \bS\ra Z$ is the contraction of the $(-2)$-curve $\bC_9$ corresponding to $C_9\subset S$, resulting
in the rational double point $b\in Z$ of type $A_1$. We have
$\Spec(k)=\bD_0\cap \bC_9=\bD_0\cap h^{-1}(b)$,
the latter by \cite[Theorem 4]{Artin1966}. According to the   Nakayama Lemma,  the induced morphism
$h\colon \bD_0\ra D_0$ is an isomorphism. In turn, $D_\red=D_0$ is the rational cuspidal curve.

The Adjunction Formula for $D\subset Z$ yields $\deg(\omega_D)=(K_Z+D)\cdot D=2$,
thus $\chi(\O_D)=-1$.
The normal surface $Z$ satisfies Serre's Condition $(S_2)$, so the Cartier divisor $D$ 
satisfies $(S_1)$. Consequently, the  ideal sheaf $\shN\subset\O_D$ for the closed subscheme $D_0\subset D$ is torsion-free.
It is invertible at $a\in D_0$, where $D_0$ is Cartier, and has rank one, hence $\shN$ is invertible as sheaf on $D_0$.
Thus $D$ is a ribbon on the rational cuspidal curve $D_0$.
The short exact sequence $0\ra\shN\ra\O_D\ra\O_{D_0}\ra 0$ yields  $\chi(\shN)=\chi(\O_D)-\chi(\O_{D_0})=-1$.
In turn,    $\deg(\shN)=-1$.
\qed

\medskip
Now consider the preimage $D'=\nu^{-1}(D)$ of the conductrix $D\subset Z$ on the normal del Pezzo surface 
with respect to the double covering $\nu\colon Z'\ra Z$.
 
\begin{proposition}
\mylabel{ribbon upstairs}
The scheme $D'$ is a ribbon on  $D'_\red=\PP^1$, for the invertible  sheaf $\shM=\O_{\PP^1}(-2)$.
The morphism $\nu\colon D'_\red\ra D_\red$ factors as
the normalization map $\PP^1\ra D_\red$ followed by the relative Frobenius   $F\colon \PP^1\ra \PP^1$.
Moreover, the induced map
\begin{equation}
\label{infinitesimal gluing}
\O_{\PP^1}(-1)^{\otimes 2}=F^*(\O_{\PP^1}(-1)) = \nu^*(\shN)\lra \shM = \O_{\PP^1}(-2)
\end{equation}
is bijective. 
\end{proposition}

\proof
Let $\zeta\in D$ be the generic point. Then  $\O_{D,\zeta}=F[\epsilon]$, where $F=\kappa(\zeta)$ is the function field of $D_\red$
and $\epsilon$ is an indeterminate subject to  $\epsilon^2=0$.  
 The induced extension $F\subset\O_{D',\zeta}/(\epsilon)$ has degree two.
Since $Z'$ is normal, the local ring $\O_{Z',\zeta}$ is a discrete valuation ring,
and the fiber $\nu^{-1}(\zeta)$ has embedding dimension at most one. It follows that the local Artin ring $\O_{D',\zeta}/(\epsilon)$
is a field, which must be purely inseparable over $F$. This shows that $D'=2D'_\red$.

The short exact sequence $0\ra\O_{Z'}(-D'_\red)\ra\O_{Z'}\ra\O_{D'_\red}\ra 0$ yields an exact sequence
$$
H^1(Z',\O_{Z'})\lra H^1(D'_\red,\O_{D'_\red})\lra H^2(Z',\O_{Z'}(-D'_\red)).
$$
According to Theorem \ref{del pezzo}, the term on the left vanishes, whereas the term on the right is Serre dual
to $H^0(Z',\O_{Z'}(-D'_\red))$, which vanishes as well. Thus $h^1(\O_{D'_\red})=0$, and it follows that
$D'_\red=\PP^1$.

The Adjunction Formula gives $\deg(\omega_{D'})=(D'-D')\cdot D'=0$, and this implies $\chi(\O_{D'})=0$.
The ideal sheaf $\shM$ for the inclusion $D'_\red\subset D$ has $\chi(\shM)=\chi(\O_{D'})-\chi(\O_{D'_\red})= -1$.
Furthermore, it is torsion-free and of rank one as sheaf on $D'_\red=\PP^1$, thus $\shM=\O_{\PP^1}(-2)$.
In turn, $D'=\nu^{-1}(D)$ is a ribbon on $\PP^1$ with respect to the dualizing sheaf $\omega_{\PP^1}=\O_{\PP^1}(-2)$.
The morphism $\nu\colon Z'\ra Z$ induces the map \eqref{infinitesimal gluing}, which must be injective.
It is bijective, because  both sides have the same Euler characteristic.
\qed

\medskip
Note that a similar situation already occurred  in the study of Beauville's Kummer varieties  
in characteristic two (\cite[Proposition 7.3]{Schroeer2009b}).
We now   clarify the flatness properties of the normal del Pezzo surface over the normal Enriques surface:

\begin{corollary}
\mylabel{flatness}
The double covering $\nu\colon Z'\ra Z$   is flat precisely over the complement of the $E_8$-singularity $a\in Z$.
\end{corollary}

\proof
According to \cite[Chapter III, Theorem 9.9]{Hartshorne1977}, our finite morphism is flat at all points
$z\in Z$ where the fiber $\nu^{-1}(z)$ has length two.
Since $Z'$ is Cohen--Macaulay, this  automatically holds when the local ring $\O_{Z,z}$ is regular
(\cite[proposition~6.1.5]{EGAIVb}), that is, for $z\neq a,b$.
In light of our description of the induced map $\nu\colon D'_\red\ra D_\red$ as  a composition of relative Frobenius
with the normalization, this fiber  over the $E_8$-singularity $a\in Z$ has length four,
whereas the fiber of the $A_1$-singularity $b\in Z$ has length two.
\qed
 
\medskip
It remains to determine the singularities on $Z'$, and then to compute the group $\APic(Z')$.
Recall that the normal Enriques surface $Z$ contains a rational double point $b\in Z$ of type $A_1$.
Let $b'\in Z'$ be the corresponding point on the normal del Pezzo surface.
Around this point, the double covering $\nu\colon Z'\ra Z$ is flat, so the
local ring $\O_{Z',b'}$ must be singular. Moreover, all singularities must be rational double points,
according to Theorem \ref{del pezzo}. We now make a preliminary observation:

\begin{lemma}
\mylabel{singularity preliminaries}
The point $b'\in Z'$ is the only singularity on the normal del Pezzo surface $Z'$ that lies on the preimage $\nu^{-1}(D)$
of the conductrix. If the Enriques surface $S$ is supersingular,  we actually have $\Sing(Z')=\{b'\}$.
\end{lemma}

\proof
The reduction    $D_0\subset Z$ of the conductrix is Cartier away from the singularity $b\in Z$, because the other
singularity   is factorial. In turn, its preimage $\nu^{-1}(D_0)\subset Z'$ is Cartier away from $b'\in Z'$.
According to Proposition \ref{ribbon upstairs}, this preimage is regular, hence the local rings $\O_{Z',x}$ are regular when $x\in Z'$
maps to $D\smallsetminus\{b\}$.

For the second statement, consider the unique genus-one fibration $\varphi\colon S\ra\PP^1$, and the normalization $S'\ra S$
of the K3-like covering. The arguments from \cite[Proposition 8.1]{Schroeer2017}  show  that
the Stein factorization of the composite map $S'\ra\PP^1$ is given by the relative Frobenius map $F\colon \PP^1\ra\PP^1$.
In particular, $S'$ is regular over   the relative smooth locus $\Reg(S/\PP^1)$.
Now suppose that $S$ is supersingular. Then there is only one multiple fiber, 
whose reduction is $C_0+\ldots+C_8$, with Kodaira symbol $\II^*$. All other geometric fibers are rational cuspidal curves,
and $C_9\subset S$ is the curve of cusps.
Since $C_1+\ldots+C_8$ and $C_9$ are contracted by $f\colon S\ra Z$, we see that the singular locus of $Z'$ lies
over $D_0=f(C_0)$.
\qed

\medskip
We now can unravel the   picture completely:

\begin{proposition}
\mylabel{singularity d5}
The rational double point $\O_{Z',b'}$ has type $D_5$, and this is the only
singularity on the normal del Pezzo surface $Z'$.
\end{proposition}

\proof
Let $r\colon X\ra Z'$ be the minimal resolution of singularities. Then $X$ is a weak del Pezzo surface of degree $K_X^2=K_{Z'}^2=4$.
It is obtained from $\PP^2$ by blowing-up $5=9-4$ points. The Picard number is $\rho(X)=6=1+5$, so the morphism $r\colon X\ra Z'$ contracts
five $(-2)$-curves $E_1,\ldots,E_5\subset X$. 
The possible configurations of rational double points  on normal del Pezzo surfaces of degree four were classified by Dolgachev
\cite[Section 8.3.6]{Dolgachev2012}. There are only two possibilities with five exceptional divisors, namely $D_5$ or
the configuration $A_3+A_1+A_1$.

Seeking a contradiction, we assume that $\Sing(Z')$ is given by $A_3+A_1+A_1$.
In light of Lemma~\ref{singularity preliminaries}, the Enriques surface $S$ must be classical.
The unique genus-one fibration $\varphi\colon S\ra\PP^1$ thus has two multiple fibers.
Without restriction, these occur over the points $0,\infty\in\PP^1$,   with respective Kodaira symbols $\II^*$ and $\II$.
Write $C_\infty$ for the reduced part of the multiple  fiber $\varphi^{-1}(\infty)$,
and   $D_\infty=f(C_\infty)$ for its image on the normal Enrique surface $Z$.
This is the rational cuspidal curve with self-intersection $D_\infty^2=1/2$.
Its    preimage on the canonical covering $\tilde{Z}\ra Z$, which coincides with the canonical covering
for the inclusion $\alpha_2\subset\Pic_{D_\infty/k}$, must be a ribbon on the projective line,
according to  \cite[Lemma 4.2]{Schroeer2017}.
In turn, the preimage  on the normal del Pezzo surface is of the form $\nu^{-1}(D_\infty)=2\Theta$, where $\Theta=\PP^1$
and rational selfintersection
$\Theta^2=1/4$.

Moreover, the Weil divisor  $\nu^{-1}(D_0)=\PP^1$ is linearly equivalent to $\nu^{-1}(D_\infty)=2\Theta$.
Since both pass through the singular point $b'\in Z'$, the local class group $\Cl(\O_{Z',b'})$ is not annihilated by two,
and we conclude that $b'\in Z'$ has type $A_3$ rather than $A_1$. Furthermore,  the germ $\Theta_{b'}$ generates the local class group.
The   five exceptional curves $E_i$ and the strict transform $\Theta^*\subset X$ of $\Theta\subset Z'$
have normal crossings, because the scheme $\Theta$ is regular. After reordering, their dual graph takes this  form:
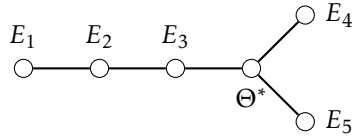
\begin{equation*}
\begin{gathered}
\begin{tikzpicture}
[node distance=1cm, font=\small]
\tikzstyle{vertex}=[circle, draw, inner sep=0mm, minimum size=1.5ex]
\node[vertex]	(E1)  	at (0,0) 	[label=above:{$E_1$}] 			{};
\node[vertex]	(E2)			[right of=E1, label=above:{$E_2$}]	{};
\node[vertex]	(E3)			[right of=E2, label=above:{$E_3$}]	{};
\node[vertex]	(Th)			[right of=E3, label=below:{$\Theta^*$}]	{};
\node[vertex]	(E4)			[above right of=Th, label=right:{$E_4$}]	{};
\node[vertex]	(E5)			[below right of=Th, label=right:{$E_5$}]	{};
 
\draw [thick] (E1)--(E2)--(E3)--(Th)--(E4);
\draw [thick] (Th)--(E5);
\end{tikzpicture}
\end{gathered}
\end{equation*}
It follows that the rational pull-back in the sense of Mumford \cite{Mumford1961} is given by 
$$
r^*(\Theta) = \Theta^* + \frac{1}{4}(   E_1+ 2E_2 + 3E_3 + 2 E_4 + 2E_5),
$$
which yields $1/4=\Theta^2= r^*(\Theta)^2 = (\Theta^*\cdot\Theta^*) + 3/4+2/4+2/4$.
Consequently $(\Theta^*\cdot\Theta^*)=-6/4$, contradicting that this selfintersection number on
the regular surface $X$ is an integer.
\qed

\begin{proposition}
\mylabel{Pic Z'}
The normal del Pezzo surface has $\APic(Z')=\ZZ$, and the subgroup $\Pic(Z)=\ZZ K_{Z'}$ has index four. In particular, the generator
$\Theta\in\APic(Z')$ has selfintersection $\Theta^2=1/4$.
\end{proposition}

\proof 
Let $r\colon X\ra Z'$ be the minimal resolution of singularities, and $E_1,\ldots,E_5\subset X$
the five exceptional curves over the rational double point $b'\in Z'$ of type $D_5$, 
and set $\Div_E(X) =\bigoplus\ZZ E_i$. Consider the commutative diagram
$$
\begin{CD}
0	@>>>	\Div_E(X)	@>>> 	\Pic(Z')\oplus\Div_E(X)	@>\pr>>	\Pic(Z')	@>>> 0\\
@.		@VV\id  V			@VV(r^*,\can) V				@VV\can V\\
0	@>>>	\Div_E(X)	@>>>	\Pic(X)			@>>r_*>	\APic(Z')	@>>> 0\\
\end{CD}
$$
The vertical maps are injective, and the one in the middle has index four, because the intersection forms
on $\Pic(X)$ and $\Pic(Z')\oplus\Div_E(X)$ have discriminants $\delta=1$ and $\delta=4^2$,
respectively. Applying the Snake Lemma, we see that the cokernel for $\Pic(Z')\subset\APic(Z')$ 
has order four. It also sits inside $\Cl(\O_{Z',b'})$, which   is cyclic
of order four. 

Thus the group $\APic(Z')$ is an extension of $\ZZ/4\ZZ$ by $\ZZ$.
It remains to check that it is torsion-free,
in other words the inclusion of $L=\Div_E(X)$ inside $\Pic(X)$ is primitive.
Consider the dual lattice $L^*=\Hom(L,\ZZ)$ and the resulting 
discriminant group $L^*/L$, which comes with a perfect $\QQ/\ZZ$-valued pairing.
The over-lattices $L\subset L'$ correspond  to  totally isotropic subgroup 
$T\subset L^*/L$, via $T=L'/L$, according to \cite[Section 4]{Nikulin1980}. In our case, the discriminant group  $L^*/L=\ZZ/4\ZZ$ is cyclic,
so there are no such subgroups. In turn, $L\subset\Pic(X)$ is primitive,
thus $\APic(Z')=\ZZ$. The generator $\Theta$ satisfies $4\Theta=K_{Z'}$, and the intersection number $\Theta^2=1/4$ follows.
\qed

\medskip
By Artin's classification \cite{Artin1977} of rational double points in positive characteristics,
there are actually two isomorphism classes   of type $D_5$, which are 
denoted by  $D_5^0$ and $D_5^1$. The former is simply-connected, the latter not.
Since $\nu\colon Z'\ra Z$ is a finite universal homeomorphism and 
rational double points of type $A_1$ are simply-connected, we see that our singularity
$\O_{Z',b'}$ is also simply-connected, whence formally given by the equation
$z^2+y^2z+x^2y=0$. 

%

\section{Cones and Fano varieties}
\mylabel{Cones}

In this section we collect some facts on cones, which complement the discussions by Grothen\-dieck \cite[\S 8]{EGAII} and  Koll\'ar \cite[Section 3.1]{Kollar2013}.
They will be used to construct Fano threefolds with unusual torsion  in the next section.

Let $k$ be a ground field of arbitrary characteristic $p\geq 0$
and $B$ be a proper connected scheme. Suppose that  $\shE$ is a locally free sheaf of rank two, and consider its projectivization
$$
P=\PP(\shE)=\Proj(\Sym^\bullet\shE).
$$
Let $f\colon P\ra B$ the structure morphism, whose fibers are copies of the  projective line $\PP^1$.
From the Formal Function Theorem, one gets a split exact sequence
\begin{equation}
\label{picard sequence}
0\lra \Pic(B)\stackrel{f^*}{\lra}\Pic(P)\stackrel{\deg}{\lra}\ZZ\lra 0,
\end{equation}
where the degree is taken fiber-wise, and the splitting is given by the tautological sheaf $\O_P(1)$.
The sections $D\subset P$ correspond to invertible quotients $\shL=\shE/\shN$,
via $D=\PP(\shL)$ and $\shL=f_*(\O_P(1)|D)$. Each section is an  effective Cartier divisor.
Moreover, to simplify the notation, 
for any line bundle $\shL$ on $B$ and any section $D$ in $P$, 
we denote the restriction $f^*(\shL)|D$  by the same symbol $\shL$.

\begin{lemma}
\mylabel{formula restrictions}
With the above notation, $\O_P(D)\simeq f^*(\shN^{\otimes-1})\otimes\O_P(1)$. In particular,
$\O_D(D)=\shN^{\otimes-1}\otimes\shL$.
\end{lemma}

\proof
Set $\shF=\O_P(1)\otimes\O_P(-D)$. Tensoring $0\ra\O_P(-D)\ra\O_P\ra\O_D\ra 0$ with $\O_P(1)$ and taking direct images gives
$0\ra f_*(\shF)\ra \shE\ra \shL\ra 0$.
The map on the right is the   quotient map $\shE\ra\shL$ defining the section $D\subset P$, hence
$f_*(\shF)=\shN$. The invertible sheaf $\shF$ has degree zero on each fiber of $f\colon P\ra P$.
Hence $\shF=f^*(\shN')$ for some invertible sheaf $\shN'$ on $B$, by the exact sequence \eqref{picard sequence}.
Finally, the Projection Formula gives $\shN'=f_*f^*(\shN')=\shN$.
\qed

\medskip
From now on, we assume that $\shE=\O_B\oplus\shL$ for some ample invertible sheaf $\shL$.
Then $P=\PP(\shE)$ contains two canonical sections, namely $A=\PP(\O_X)$ and $E=\PP(\shL)$, coming from the two 
projections $\shE\ra\O_X$ and $\shE\ra\shL$, and we have
\begin{equation}
\label{canonical restrictions}
\O_A(A)=\shL\quadand \O_E(E)=\shL^{\otimes-1}.
\end{equation}
Therefore we say that $E\subset S$ is the \emph{negative section} and  $A\subset P$ is the \emph{positive section}.
Clearly, the stable base locus for the invertible sheaf $\O_P(A)$ is contained in $A$,
and the restriction $\O_A(A)$ is ample. According to Fujita's result 
(see \cite{Fujita1983}, compare also \cite{Ein2000}),
the invertible sheaf $\O_P(A)$ is semiample, and we get a contraction
$$
r\colon P\lra X=\Proj \bigoplus_{i\geq 0} H^0(P,\O_P(iA))
$$
to some projective scheme $X$, with  $\O_P(nA)=r^*\O_X(n)$
for some $n\geq 0$ sufficiently large, and $\O_X=r_*(\O_P)$. Clearly, the connected closed set $E\subset P$ is mapped
to a closed point $x_0\in X$. The \emph{exceptional set} for the morphism $r\colon P\ra X$
is defined as the closed set $\Exc(P/X)=\Supp(\Omega^1_{P/X})$.

\begin{lemma}
\mylabel{exceptional set}
The exceptional set $\Exc(P/X)$ for the morphism $r\colon P\ra X$ coincides with the negative section $E=\PP(\shL)$.
\end{lemma}

\proof
The exceptional set is the union of all irreducible curves $C\subset P$ that are disjoint from the positive section 
$A=\PP(\O_B)$. We have to check that each such $C$ is contained in the negative section $E=\PP(\shL)$.
For this, we may pass to the base-change $C'\times_BP$, where $C'\ra C$ is the normalization,
and assume that $B$ is an irreducible regular curve and $C\subset P$
is a section. Then $P$ is a regular surface, and both $E,C\subset P$ are mapped to points in $X$.
Recall that $E\subset P$ corresponds to the invertible quotient $\shL=\shE/\shN$, where $\shN=\O_B$.
By the Hodge Index Theorem, we have $\deg(\shL)-\deg(\shN)=E^2<0$, and the same holds for $C\subset P$. 
Since such an invertible quotient
$\shL=\shE/\shN$ is unique, $C=E$ follows.
\qed

\medskip
In turn, the morphism $r\colon P\smallsetminus E\ra X\smallsetminus\{x_0\}$ is an isomorphism.
One also says that $X$ is the \emph{projective cone} on $B$ with respect to the ample sheaf $\shL$, with
vertex $x_0\in X$.
By abuse of notation, we regard the positive section $A=\PP(\O_B)$ of the $\PP^1$-bundles $P=\PP(\shE)$
also as  an ample Cartier divisor $A\subset X$.

\begin{proposition}
\mylabel{pic cone}
The projective scheme $X$ has Picard scheme $\Pic_{X/k}=\ZZ$, and this is generated
by the ample sheaf  $\O_X(A)$.
\end{proposition}

\proof
Let $\shF$ be an invertible sheaf on $X$. Its preimage takes the form $r^*(\shF)=f^*(\shN)\otimes\O_P(dA)$
for some unique $\shN\in\Pic(B)$ and $d\in\ZZ$. Since $r^*(\shF)|E=\O_E$ and $\O_E(dA)=\O_E$,
we have $\shN=\O_D$, and hence $\Pic(X)=\ZZ$, generated by $\O_P(A)$.
The statement on the Picard scheme follows likewise, by working with the infinitesimal extension
$B\otimes_kk[\epsilon]$.
\qed

\begin{proposition}
\mylabel{canonical cone}
If $B$ is Gorenstein, then  $P$ is Gorenstein, and  we  have the Canonical Bundle Formula
$$\omega_P=\O_P(-2E)\otimes f^*(\omega_B\otimes\shL^{\otimes -1}).$$
\end{proposition}

\proof
The scheme $P$ must be Gorenstein by \cite{Watanabe-Ishikawa-Tachibana-Otsuka1969}.
In light of the exact sequence \eqref{picard sequence}, the dualizing sheaf has the form  
$\omega_P=\O_P(-2E)\otimes f^*(\omega_B\otimes\shN)$ for some
invertible sheaf $\shN$ on $B$.  Since $\O_E(E)=\shL^{\otimes -1}$,
the Adjunction Formula for the effective Cartier divisor $E\subset P$ yields
$$
\omega_B=(\omega_P\otimes\O_P(E))|E= \shL^{\otimes 2}\otimes\omega_B\otimes\shN\otimes\shL^{\otimes-1}.
$$
The assertion follows.
\qed
 
\medskip
Recall that a scheme $V$ is called a \emph{Fano variety} if it is proper, with $h^0(\O_V)=1$,
all local rings $\O_{V,a}$ are Gorenstein, and the dualizing sheaf $\omega_V$ is
anti-ample. Note that we make no other assumptions on the singularities.
Fano varieties of dimension $n\geq 1$ come  with two important numerical invariants:
the \emph{degree} and the \emph{index}
$$
\deg(V)=(-K_V)^n=c_1^n(\omega_V^{\otimes-1})>0\quadand \ind(V)>0.
$$
The latter is defined as  the divisibility of  $\omega_V$  in the numerical  group $\Num(V)$.

\begin{theorem}
\mylabel{cone fano}
Suppose the following:
\begin{enumerate}
\item The projective scheme $B$ is a Fano variety of dimension $n\geq 1$.
\item The ample sheaf $\shL$ has the property $\shL^{\otimes m}=\omega_B^{\otimes -1}$ for some $m\geq 1$.
\item The group $H^i(B,\shL^{\otimes t})$ vanishes for all integers $i,t\geq 1$.
\end{enumerate}
Then  the projective scheme $X$ is Gorenstein, and we have the 
the  Canonical Bundle Formula
$K_X=-(m+1)A$. In particular, $X$ is a Fano variety of dimension $n+1$, with 
numerical invariants
$$
\deg(X) = \frac{(m+1)^{n+1}}{m^n} \deg(B)\quadand \ind(X) = m+1.
$$
\end{theorem}

\proof
Assumption (iii) guarantees that the cone $X$ is Cohen--Macaulay, according to  \cite[Proposition~3.13]{Kollar2013}.
In light of Condition (ii) and Proposition \ref{pic cone}, the dualizing sheaf on the $\PP^1$-bundle
is
$$\omega_P=\O_P(-2E)\otimes f^*(\shL^{\otimes -m-1}).$$
Furthermore, we have $f^*(\shL)=\O_P(A-E)$, which follows from \eqref{canonical restrictions} 
and the exact sequence  \eqref{picard sequence}. This gives $K_P= -(m+3)E - (m+1)A$.
Since $A$ is disjoint from the exceptional locus $E$, we see that $X$ is Gorenstein, having 
$K_X= -(m+1)A$. With Proposition \ref{pic cone} we conclude that $X$ is a Fano variety
with index $m+1$. The degree is  
$(-K_X)^{n+1} = (m+1)^{n+1} A^{n+1}$.
In light of \eqref{canonical restrictions}, we have $A^{n+1}= c_1^n(\shL) = \frac{1}{m^n}(-K_B)^n$.
This concludes the proof.
\qed

\section{Fano threefolds with unusual torsion}
\mylabel{Fano threefolds}

Let $k$ be an algebraically closed ground field of characteristic $p=2$,
and $S$ be a simply-connected Enriques surface endowed with an ADE-curve $E\subset X$
as in Theorem \ref{del pezzo}. This actually exists, as we saw in Section \ref{Exceptional enriques}.
The resulting contraction $f:S\ra Z$ yields a normal Enriques surface,
coming with a $G$-torsor $\tZ\ra Z$, with respect to the Cartier dual $G$ for the unipotent group scheme $P=\Pic^\tau_{Z/k}$ of
order two. As explained in Section \ref{Exceptional enriques},
the normalization $Z'\ra\tZ$  gives a del Pezzo surface 
with $\Pic(Z')=\ZZ K_{Z'}$ and $h^1(\O_{Z'})=0$, of  degree $K_{Z'}^2\in\{2,4,6,8\}$.
It comes with an inseparable  double covering   $\nu:Z'\ra Z$.

We now consider the    $\PP^1$-bundle $P=\PP(\shE)$ with $\shE=\O_{Z'}\oplus\omega_{Z'}^{\otimes-1}$,
and the ensuing contraction 
$$
r:P\lra X=\Proj\bigoplus_{i\geq 0} H^0(P,\O_P(iA))
$$
of the negative section $E=\PP(\omega_{Z'}^{\otimes-1})$, defined via
the positive section $A=\PP(\O_{Z'})$, as explained in the previous section.
The singular loci  can be written as
$$
\Sing(P)=\PP^1_{\Sing(Z')}\quadand \Sing(X)=r(\PP^1_{\Sing(Z')}).
$$
 
\begin{proposition}
\mylabel{normal threefold}
The scheme $X$ is a normal Fano threefold of degree $(-K_X)^3=8\cdot K_{Z'}^2$, index $\ind(X)=2$,
and $h^i(\O_X)=0$ for all integers $i\geq 1$.
The Picard scheme is $\Pic_{X/k}=\ZZ$, which is generated by $\O_X(A)$,  
such that  $\omega_X=\O_X(-2A)$.
\end{proposition}

\proof
The statement about the Picard scheme follows from Proposition \ref{pic cone}.
According to Corollary \ref{raynaud vanishing}, we have $h^1(\omega_{Z'}^{\otimes t})=0$.
We now can apply Theorem \ref{cone fano} with $B=Z'$,  $\shL=\omega_{Z'}$ and $m=1$
and conclude that $X$ is a Fano threefold of index two and degree $(-K_X)^3=8\cdot K_{Z'}^2$.

It remains to compute the cohomological invariants. The group $H^1(X,\O_X)$ 
is the Lie algebra for the Picard scheme, hence vanishes.
Since $X$ is integral and $\omega_X^{\otimes-1}$
is ample, the group $H^0(X,\omega_X)$ vanishes. Serre duality gives $h^3(\O_X)=0$.
To see that $h^2(\O_X)$ vanishes, it suffices to check $\chi(\O_X)=1$. The Leray--Serre spectral sequences for the 
fibration $f:P\ra Z'$ and the contraction $r:P\ra X$ give
$$
1= \chi(\O_P) = \chi(\O_X) - \chi( R^1r_*(\O_P)) + \chi(R^2r_*(\O_P)).
$$
Recall that $E\subset P$ is a negative section that is contracted.
The short exact sequence
$$0\ra\O_E(-nE)\ra\O_{(n+1)E}\ra\O_{nE}\ra 0$$
yields
$$
H^i(\O_E(-nE))\lra H^i(\O_{(n+1)E})\lra H^i(\O_{nE})\lra H^{i+1}(\O_E(-nE)).
$$
Using the identification $E=Z'$ and $\O_E(-E)=\omega_{Z'}^{\otimes-1}$, we see that $H^2(\O_E(-nE))$ 
vanishes for all $n\geq 1$.
According to Corollary~\ref{raynaud vanishing}, the groups $H^1(\O_E(-nE))$ vanish as well.
With the Formal Function Theorem  we conclude that $R^1r_*(\O_P)= R^2r_*(\O_P)=0$,
and hence $\chi(\O_X)=1$.
\qed

\medskip
Using the   inclusion $Z'\subset X$ coming from $Z'\subset P\ra X$
and the   inseparable double covering $\nu:Z'\ra Z$,
we now  form the   cocartesian and cartesian square
\begin{equation}
\label{conductor square}
\begin{CD}
Z'		@>\can>>	X\\
@V\nu VV		@VVV\\
Z		@>>>	Y.
\end{CD}
\end{equation}
Then  $Y=X\amalg_ZZ'$ is an integral proper  threefold, with normalization $X$.
A priori, the amalgamated sum exists
as an algebraic space (\cite[Theorem 6.1]{Artin1970}). Since the normalization 
map $\nu:X\ra Y$ is a universal homeomorphism, the algebraic space $Y$ must  be a scheme 
(\cite[Theorem 6.2.2]{Olsson2016}). It contains the normal Enriques surface $Z$ as
a closed subscheme.
By construction, the singular locus  is given by $\Sing(Y)=Z\cup r(\PP^1_{\Sing(Z')})$.   

\begin{proposition}
\mylabel{good threefold}
The integral proper threefold $Y$ has the following properties:
\begin{enumerate}
\item The dualizing sheaf $\omega_Y$ is invertible and anti-ample, with   $(-K_Y)^3=K_{Z'}^2$.
\item We have  $\Num(Y)=\ZZ$.
\item The Euler characteristic is $\chi(\O_Y)=1$.
\item For each closed point  $y\in Y$, the local ring $\O_{Y,y}$ satisfies Serre's Condition $(S_2)$.
It is  actually Cohen--Macaulay provided that   $\O_{Z,y}\subset\O_{Z',y}$ is flat.
\end{enumerate}
\end{proposition}

\proof 
The conductor square \eqref{conductor square} yields the  short exact sequence 
\begin{equation}
\label{ses normalization}
0\lra\O_Y \lra \O_X \oplus \O_Z\lra \O_{Z'} \lra 0.
\end{equation}
The normal Enriques surface $Z$ and the normal del Pezzo surface $Z'$ both have Euler characteristic $\chi=1$,
and the same holds for the normal Fano threefold $X$, by
Proposition \ref{normal threefold}. This gives $\chi(\O_Y)=1$.

Since the map $\nu:X\ra Y$ is surjective, each integral curve on $Y$ is the image of an integral curve on $X$,
hence  the induced map 
$\nu^*\colon \Num(Y)\to \Num(X)$ is injective, and it follows that the group $\Num(Y)$ is free of rank one.

The singular locus $\Sing(Y)=Z\cup r(\PP^1_{\Sing(Z')}) $ consists of   an irreducible surface and a curve.
Let $\zeta\in Y$ be the generic point of
the conductor locus $Z$, and write $\zeta'\in X$ for the corresponding point on $X$.
Then $\kappa(\zeta)\subset \kappa(\zeta')$ is an inseparable field extension of degree two,
and the subring $\O_{Y,\zeta}\subset \O_{X,\zeta}$ comprises all ring elements
whose class in the residue field $\kappa(\zeta')$ lies in the subfield $\kappa(\zeta)$.
It follows that the local ring $\O_{Y,\zeta}$ is Gorenstein, compare the discussion in \cite[Appendix A]{Fanelli-Schroeer2020}.
Thus the dualizing sheaf $\omega_Y$ is invertible on some open subset containing all 
codimension-one points.
Now let $y\in Y$ be an arbitrary point, and write $x\in X$ for the corresponding point.
Since both $\omega_X$ and $\O_X(Z')$ are invertible, 
\cite[Proposition~A.4]{Fanelli-Schroeer2020} applies, and we conclude that the dualizing sheaf $\omega_Y$ is invertible.

According to Theorem \ref{cone fano}, the dualizing sheaf on the Fano threefold $X$ is given by $\omega_X=\O_X(-2Z')$.
The relative dualizing sheaf for the finite birational morphism $\nu:X\ra Y$ is defined 
by the equality $\nu_*(\omega_{X/Y})=\uHom(\nu_*(\O_X),\O_Y)=\nu_*\O_X(-Z')$.
From $\omega_X=\nu^*(\omega_Y)\otimes\omega_{X/Y}$ we conclude $\nu^*(\omega_Y)=\O_X(-Z')$.
In particular, $K_Y^3 = (-Z')^3 = -K_{Z'}^2$ holds by Lemma \ref{formula restrictions}.

Now fix a closed point $y\in Z$, and consider the three-dimensional local ring $\O_{Y,y}$.
It is Cohen--Macaulay by Theorem \ref{cone fano}, provided that $y\not\in Z'$.
Now suppose that $y\in Z$.
The short exact sequence \eqref{ses normalization} induces a long exact sequence
\begin{equation}
\label{sequence local cohomology}
H^{i-1}_y(\O_{Z'})\lra H^i_y(\O_Y)\lra H^i_y(\O_X)\oplus H^i_y(\O_Z)\lra H^i_y(\O_{Z'}) 
\end{equation}
of local cohomology groups. The two-dimensional schemes $Z$ and $Z'$ are Cohen--Macaulay,
so their local cohomology groups vanish in degree $i<2$.
Likewise, the three-dimensional scheme $X$ is Cohen--Macaulay, so we have vanishing in degree $i<3$.
It follows that $H_y^i(\O_Y)=0$ for $i\leq 1$, hence the local ring $\O_{Y,y}$ satisfies Serre's Condition $(S_2)$.

Finally suppose that $\O_{Z,y}\subset\O_{Z',y}$ is flat.
Then the inclusion is a direct summand, in particular the
induced map $H^2_y(\O_Z)\ra H^2_y(\O_{Z'})$ is injective. This ensures $H^2_y(\O_Y)=0$,
so the local ring $\O_{Y,y}$ is Cohen--Macaulay.
\qed

\medskip
We see that the scheme $Y$ qualifies as a Fano variety, except that it is not necessarily Cohen--Macaulay.
By definition, a local noetherian ring   is called \emph{Gorenstein} if it is Cohen--Macaulay, and the dualizing
module is invertible. Without the former condition, one should use the term \emph{quasi-Gorenstein}  instead.
Thus it is natural to call our scheme $Y$ a \emph{quasi-Fano  variety}, or a \emph{Fano variety that is not necessarily Cohen--Macaulay}.
Our construction yields unusual torsion in the Picard scheme:

\begin{theorem}
\mylabel{exotic threefold}
The  integral  Fano threefold $Y$  that is not necessarily Cohen--Macaulay has the following property:
$$
\Upsilon_{Y/k}=\Pic^\tau_{Y/k}=
\begin{cases}
\ZZ/2\ZZ	& \text{if the Enriques surface $S$ is classical;}\\
\alpha_2	& \text{if $S$ is supersingular.}
\end{cases}
$$
Moreover,  $h^1(\O_Y)=h^2(\O_Y)=0$ in the former case, and
$h^1(\O_Y)=h^2(\O_Y)=1$ in the latter.
\end{theorem}

\proof
The conductor square \eqref{conductor square} yields a short exact sequence of multiplicative abelian sheaves
$$1\ra\O_Y^\times\ra\O_X^\times\oplus\O_{Z}^\times \ra\O_{Z'}^\times\ra 1$$
as explained in \cite[Proposition 4.1]{Schroeer-Siebert2002}.
This holds not only in the Zariski topology, but also in the finite flat topology.
In turn, we get an exact sequence of group schemes
$$
0\lra \Pic_{Y/k}\lra\Pic_{X/k}\oplus\Pic_{Z/k}\lra \Pic_{Z'/k}.
$$
The map on the left is indeed injective, because $H^0(Z,\O_Z)\ra H^0(Z',\O_{Z'})=k$ is surjective.
Since $\rho(Y)=1$, an invertible sheaf $\shL$ on $Y$ is numerically trivial if and only
if its restriction to $X$, or equivalently to $Z$, is numerically trivial. Since $\rho(X)=1$,
the same holds for invertible sheaves on $X$ and their restriction to $Z'$. In turn, we get an 
induced exact sequence 
$$
0\lra \Pic^\tau_{Y/k}\lra\Pic^\tau_{X/k}\oplus\Pic^\tau_{Z/k}\lra \Pic^\tau_{Z'/k}.
$$
The term for the   del Pezzo surface $Z'$ vanishes, by Theorem \ref{del pezzo}.
Also the term for the normal Fano threefold $X$ is zero, according to Proposition  \ref{normal threefold}.
We thus get an identification $\Pic^\tau_{Y/k}=\Pic^\tau_{Z/k}$.
We saw in Proposition \ref{numerical group} that the morphism $S\ra Z$ from the simply-connected
Enriques surface $S$ to the normal Enriques surface $Z$ induces an identification $\Pic^\tau_{Z/k}=\Pic^\tau_{S/k}$
and the assertion on Picard schemes and its maximal unipotent quotient follows.

It remains to check the assertions on $h^i(\O_Y)$. This follows from the long exact sequence
for \eqref{ses normalization}, together with  $h^i(\O_X)=h^i(\O_{Z'})=0$ for $i\geq 1$. 
\qed

\medskip
Let us now examine a concrete example for the construction of $Y$:
Suppose that our simply-connected Enriques surface $S$ is an exceptional Enriques surface of type $T_{2,3,7}$,
and that $f:S\ra Z$ is the contraction of the 
ADE-curves $C_1+\ldots+C_8$ and $C_9$, as analyzed in Section \ref{Exceptional enriques}.
Recall that $\Sing(Z)=\{a,b\}$ and $\Sing(Z')=\{b'\}$, 
where  the corresponding local rings are rational double points of type $E_8, A_1$ and $D_5$, respectively.
Moreover, the resulting threefold $Y$ has $\Sing(Y) = Z\cup r(\PP^1_{b'})$.
It turns out that the closed point $a\in Y$, which comes from the $E_8$-singularity $\O_{Z,a}$, plays a special role:
 
\begin{proposition}
In the above situation, our     quasi-Fano threefold $Y$ has degree $(-K_Y)^3=4$ and  $\Num(Y)=\ZZ K_Y$.
Moreover,   $Y$ is Cohen--Macaulay outside $a\in Y$, whereas the local ring $\O_{Y,a}$ 
satisfies $(S_2)$ but not $(S_3)$.
\end{proposition}

\proof
By Corollary \ref{flatness}, the double covering $Z'\ra Z$ is flat precisely outside $a\in Z$. So Proposition \ref{good threefold} tells us that 
the local rings $\O_{Y,y}$ are Cohen--Macaulay for $y\neq a$.
It remains to understand the case $y=a$.
Now the cokernel $M$ for the inclusion $\O_{Z,a}\subset\O_{Z',a}$ still is torsion-free of rank one, but
fails to be  invertible. Since the local ring $\O_{Z,a}$ is a factorial, the bidual $M^{\vee\vee}$ is invertible,
and the cokernel $F=M^{\vee\vee}/M$ is finite and non-zero. Using the long exact sequence
for the short exact sequence $0\ra M\ra M^{\vee\vee}\ra F\ra 0$, one  easily infers $H^1_a(M)\neq 0$.
With the exact sequence \eqref{sequence local cohomology} we infer that the map $H^2_a(\O_Z)\ra H^2_a(\O_{Z'})$ is not injective, and hence $H^2_a(\O_Y)\neq 0$.
In turn, the local ring $\O_{Y,a}$ is not Cohen--Macaulay.

Combining Proposition \ref{picard group downstairs} and \ref{good threefold} we get the degree $(-K_Y)^3= K_{Z'}^2=4$. This is not a multiple of eight,
and it follows that $\Num(Y)$ is generated by the canonical class. 
\qed

\medskip
By construction, our integral Fano threefold $Y$ is not Cohen--Macaulay in codimension three,
and not regular in codimension one.
In light of Salomonsson's equations \cite{Salomonsson2003}, the
construction works over any ground field   of characteristic $p=2$.
It would be interesting to know if there are imperfect   fields $F$ over
which $Y$ admits a twisted form $Y'$ whose local rings are normal, $\QQ$-factorial klt singularities.
Such twisted forms could appear as generic fibers in Mori fiber spaces. 
Recall that there are indeed examples of three-dimensional normal, $\QQ$-factorial  terminal singularities that are  
 not Cohen--Macaulay \cite{Totaro2019}.
In \cite{Schroeer2007}, related  problems where considered for   non-normal del Pezzo surface.
 Del Pezzo surfaces over $F$ with $\pdeg(F)=1$ were studied in \cite{Fanelli-Schroeer2020}.
Here the \emph{p-degree} is defined as $\pdeg(F)=\dim_F(\Omega^1_{F/F^p})$.


\ifx\undefined\bysame
\newcommand{\bysame}{\leavevmode\hbox to3em{\hrulefill}\,}
\fi

\end{document}